\documentclass[10pt,a4paper]{article}
\usepackage{t1enc}
\usepackage[latin1]{inputenc}
\usepackage[german,english]{babel}
\usepackage{amsfonts}
\usepackage[all]{xy}
\usepackage{algorithm2e}
\newtheorem{defi}{Definition}[section]
\newtheorem{prop}[defi]{Proposition}
\newtheorem{theo}[defi]{Theorem}
\newtheorem{lem}[defi]{Lemma}
\newtheorem{coro}[defi]{Corollary}

\begin{document}

%%references
\newcommand{\Ar}{MR0019087}
\newcommand{\Be}{MR2032983}
\newcommand{\Bi}{MR1815219}
\newcommand{\Kr}{MR1888796}
\newcommand{\BKLa}{MR1870512}
\newcommand{\BKLb}{MR1654165}
\newcommand{\BV}{MR2202556}
\newcommand{\CP}{MR0132791}
\newcommand{\CW}{MR1942303}
\newcommand{\De}{MR2787454}
\newcommand{\DD}{MR2463428}
\newcommand{\DP}{MR1710165}
\newcommand{\DKL}{MR2987387}
\newcommand{\Di}{MR2004479}
\newcommand{\EEO}{MR3030521}
\newcommand{\EM}{MR1315459}
\newcommand{\EC}{MR1161694}
\newcommand{\GL}{MR0251712}
\newcommand{\KKL}{MR1465024}
\newcommand{\Lee}{Lee10}
\newcommand{\LZ}{MR0445901}
\newcommand{\LO}{MR2340901}
\newcommand{\MP}{MR0264062}
\newcommand{\Pa}{Pi00}
\newcommand{\Pb}{MR1842395}
\newcommand{\Simon}{MR563704}
\newcommand{\Vass}{Vass}

\title{{\sc Logspace computations for Garside groups of spindle type}}
\author{Murray Elder and Arkadius Kalka\footnote{Research supported by Australian Research Council grants DP110101104  and FT110100178}}

%\subjclass[2010]{20F65, 68Q15}
%\keywords{Logspace normal form; logspace conjugacy problem; braid group; Garside group}
 %\thanks{Research supported by Australian Research Council grants DP110101104  and FT110100178}

\date{\small \today. Mathematics Subject Classification(2010): 20F36, 20F65, 68Q15. Keywords: Logspace normal form, logspace conjugacy problem, braid group, Garside group}

\maketitle

\begin{abstract}
M. Picantin introduced the notion of Garside groups {\em of spindle type}, generalizing the 3-strand braid group. % as prototypical example.
We show that, for linear Garside groups of spindle type, a normal form and a solution to the conjugacy problem are logspace computable.
For linear Garside groups of spindle type with {\em homogenous} presentation we compute a geodesic normal form in logspace.
%Furthermore, we show how to compute several other classical 3-braid normal forms in logspace, and we provide some examples of normal forms in $B_4$
%and $B_3$ that are not L-computable.
\end{abstract} 

\section{Introduction and Outline}
The word problem,  conjugacy problem and  geodesic problem are examples for classical problems in combinatorial group theory.
Here we deal with the question which of these problems can be solved in {\it logarithmic space}. Our objects of interest are braid and
Garside groups, certain generalizations of braid groups which admit a similar solution to the word and conjugacy problem. 
Since braid  groups are linear \cite{\Bi, \Kr} the word problem is solvable in logspace \cite{\LZ}.
Closely related to the word problem is the computation of a normal form. Here the problem seems to be much harder. In \cite{\EEO} it is shown that torus knot groups $\langle a,b \ | \ a^m=b^n\rangle$ for $m,n$ positive integers have logspace normal forms, which includes the 3-strand braid group. Here we prove that 
Garside normal forms are logspace computable  for a small class of simple Garside groups which includes $B_3$.
Also, we use this result to solve the conjugacy problem and the shortest word problem for some special Garside groups in logspace.
For $n\ge 4$, the braid groups are not of spindle type, and it is not known whether the Garside normal form (or any other normal form)  is logspace computable. More generally, it is an open problem whether all linear (respectively automatic) groups have logspace computable normal forms, and also whether all biautomatic groups have conjugacy problem solvable in logspace, so since braid groups are biautomatic \cite{\EC} and linear, they serve as important test cases.

In related work Diekert,  Kausch and  Lohrey consider logspace geodesic normal forms for  right-angle Artin groups and Coxeter groups \cite{\DKL}, and S. Vassileva showed that the conjugacy problem for a large class of wreath products is solvable in logspace \cite{\Vass}.
\\

%cite: Vassilva (wreath prod), Diekert (graph products, RA artin groups), EO (nilpotent cp)

\noindent{\bf Outline.} In Section \ref{II} we introduce the basic notion of an {\it logspace transducer} and recall some basic examples and 
properties of logspace computable functions and normal forms. Section \ref{III} deals with normal forms in Garside systems of spindle type. 
First we recall the definition of Garside systems (section \ref{3.1}) and their normal forms (section \ref{3.2}). In section \ref{GarsideBraids} 
we refer to the basic Garside structures in braid groups. Section \ref{3.4} introduces the main notion of so-called Garside groups {\emph of spindle type}. 
The main results on normal forms in Garside groups of spindle type are given in sections \ref{3.5} and \ref{IV}. In particular, we show that the Garside normal form is logspace computable for linear Garside groups of spindle type (section \ref{3.5.1}). Furthermore, we show that a geodesic normal form is logspace computable for linear Garside groups of spindle type with homogenous presentation (section \ref{3.5.2}). A logspace solution for the conjugacy problem in linear  Garside groups of spindle type is given in section \ref{IV}. 
In the final section we briefly consider a normal form for  $B_4$, which we show cannot  be computed in logspace. 

\section{Basic properties of logspace normal forms} \label{II}

We start with a precise definition of logspace computation, taken from \cite{\EEO}.

\begin{defi} \label{LTransduc} \hspace{-0.15cm}{\bf .}
A {\sc deterministic logspace transducer} consists of a finite state
control and three tapes: the first {\sc input tape} is read only, and
stores the input word; the second {\sc work tape} is read-write, but is
restricted to using at most $c\log n$ squares, where $n$ is the length
of the word on the input tape and $c$ is a fixed constant; and the
third {\sc output tape} is write-only, and is restricted to writing left
to right only. A transition of the machine takes as input a letter of the
input tape, a state of the finite state control, and a letter on the
work-tape. On each transition the machine can modify the work tape,
change states, move the input read-head, and write at most a fixed constant number of letters to the
output tape, moving right along the tape for each letter printed.
\end{defi}

Since the position of the read-head of the input tape is an integer between 1 and $n$, we can store it in  binary on the work tape.

\begin{defi} \label{LComput} \hspace{-0.15cm}{\bf .}
Let $X,Y$ be finite alphabets. Let $X^*$ denote the set of all finite length strings in the letters of $X$, including the empty string $\lambda$.
We call $f:X^*\rightarrow Y^*$ a
 {\sc logspace computable function}  (short: L-computable function) if there is a deterministic  logspace transducer that on input $w\in X^*$ computes $f(w)$.
\end{defi}

In the sequel we make use of the following results on L-computable functions and normal forms.
The first result is due to Lipton and Zalcstein \cite{\LZ}.
\begin{prop} \label{LZ} \hspace{-0.15cm}{\bf .}
All linear groups (groups of matrices with entries from
a field of characteristic zero) have word problem solvable in logspace.
\end{prop} Simon \cite{\Simon} extends this to linear groups over arbitrary fields.

\begin{prop} \label{freeGpNFinL} \hspace{-0.15cm}{\bf .} (see Prop. 1 in \cite{\EEO})
Let $\langle a_1, \dots, a_k \ | \ -\rangle$ be the free group of  finite rank $k$ with normal form the set of all freely reduced words over $X= \{a_i^{\pm 1}\}$. Then there is a logspace computable function $f: X^* \rightarrow X^*$ such that $f(w)$ is the normal form word for $w$.
\end{prop}
This  result  can be traced back to  \cite{\LO}. 
The next result shows that logspace computable functions are closed under composition. 
\begin{prop} \label{composition} \hspace{-0.15cm}{\bf .} (see Lemma 2 in \cite{\EEO})
If $f: X^* \rightarrow Y^*$ and $f: Y^* \rightarrow Z^*$ can both be computed in logspace,
then their composition $f \circ g: X^* \rightarrow Z^*$ can also be computed in logspace.
\end{prop}

Lastly, we note that finding a logspace computable normal form for one finite generating set implies the existence of a logspace normal form for any other  (see Prop. 5 in \cite{\EEO}).

\section{Normal forms for Garside systems of spindle type}  \label{III}

\subsection{Basic definitions for Garside systems} \label{3.1}
This introductory subsection follows in exposition and terminology closely \cite{\Pb} and \cite{Go07}. \par 
Let $M$ be a monoid. For $a, b \in M$, we say that $a$ is a {\it left divisor of $b$}, and $b$ is a {\it right multiple of $a$}, denoted by
$a\preceq b$, if there exists an element $c \in M$ such that $b=ac$. 
Similarly, for $a, b \in M$, we say that $a$ is {\it right divisor of $b$}, and $b$ is a {\it left multiple of $a$}, denoted by $b\succeq a$, if there exists an element $c \in M$ such that $b=ca$. \\
A monoid $M$ is {\it atomic} (or {\it Noetherian}) if for every element $a \in M$ there exists an integer $n$
such that the element cannot be expressed as the product of more than $n$ elements distinct from 1. An element $a\ne 1$ in $M$ is called an {\it atom} if $a=bc$ implies $b=1$ or $c=1$. An atomic monoid is necessarily generated by its atoms. Indeed, the subsets of $M$ that generate $M$ are
exactly those subsets that include the set of all atoms. 
Furthermore, if $M$ is a finitely generated atomic monoid then relations $\preceq$ and $\succeq$ are partial orders, and every element of $M$ admits only finitely many left and right divisors \cite{\DP}. We call an element $\Delta \in M$ {\it balanced} if the sets of left and right divisors coincide. In that case, this set is denoted by ${\rm Div}(\Delta )$. \\
A monoid $M$ is {\it cancellative} if for all $a, b, b', c \in M$ $abc=ab'c$ implies $b=b'$. \\
An element $m \in M$ is a {\it least common right multiple} (or a {\it right lcm})
{\it of a and b} if it is a right multiple of both $a$ and $b$, and $a \preceq c$ and $b \preceq c$ implies $m \preceq c$. \\
A {\it left lcm} is defined symmetrically using the relation $\succeq$. 
If they exist, right and left lcm's are, by definition, unique, and denoted by $a \vee b$ and $a \, \tilde{\vee }\, b$, respectively.
If $a\vee b$ exists, and $M$ is (left) cancellative, there exists a unique element c such that $a\vee b=ac$. This element is called the {\it right complement of $a$ in $b$}, and it is denoted by $a\backslash b$.
We define the {\it left complement} symmetrically. In particular, we have $a\vee b = a(a\backslash b) = b(b\backslash a)$, and $a \, \tilde{\vee }\, b= (b/a)a = (a/b)b$.

\begin{defi} \label{GaussMon} \hspace{-0.15cm}{\bf .}
A monoid $M$ is an {\sc lcm monoid} if it is Noetherian, cancellative, and every
pair of elements $a, b \in M$ admits a right and a left lcm.
\end{defi}

In an lcm monoid, for every pair of elements $(a, b)$, the set of common left divisors of $a$ and $b$ is finite and it admits a right lcm,
which is therefore the {\it greatest common left divisor} (or {\it left gcd}) {\it of $a$ and $b$}, denoted by $a \wedge b$.
The {\it right gcd} $a \, \tilde{\wedge } \, b$ is defined symmetrically. An lcm monoid is a lattice for the left and right divisibility relations
$\preceq$ and $\succeq$. \\
If $M$ is a lcm monoid, then $M$ satisfies Ore's conditions \cite{\CP}, and it embeds in
its group of right fractions, and, symmetrically, in its group of left fractions. These two groups coincide, and therefore
we may speak of the group of fractions of an lcm monoid.

\begin{defi} \label{GarsSys} \hspace{-0.15cm}{\bf .} \cite{\Pb ,Go07}
Let $G$ be a group. Denote by $G^+$ the submonoid of $G$ generated by a given subset $S$ of $G$. The pair $(G,S)$ is called a {\sc Garside system} if $G^+$ is an  lcm monoid, $G$ is its group of fractions, and there exists a balanced element $\Delta \in G^+$ such that $S={\rm Div}(\Delta )$ is finite and it generates $G$. 
Then we call $G$ a {\sc Garside group}, $G^+$ a {\sc Garside monoid}, and $\Delta $ a {\sc Garside element}. 
The elements of $S$ are called {\sc simple elements}.
\end{defi} 
$S$ is closed under $\backslash , /, \vee , \tilde{\vee }, \wedge , \tilde{\wedge }$. 
Indeed, $S$ is the closure of the atoms of $G^+$ under $\backslash $ and $\vee $. 
The functions $\partial: a \mapsto a\backslash \Delta$ and $\tilde{\partial }: a \mapsto \Delta /a$ map $G^+$ onto $S$, and the restrictions $\partial |_S$, $\tilde{\partial }|_S$ are bijections of $S$ satisfying $\tilde{\partial }|_S=(\partial |_S)^{-1}$. In particular, we have $\partial ^2(a)=\tau (a)$ and $\tilde{\partial }^2=\tau ^{-1}(a)$ for all $a\in S$, where $\tau $ denotes the inner automorphism of $G$ defined by $a\mapsto \Delta ^{-1}a\Delta $. \\
The partial orders $\preceq $ and $\succeq$ on $G^+$ naturally extend to partial orders on the whole Garside group $G$. Furthermore, we define 
for $a, b\in G$, $a \le b$, if there exist $c, c'\in G^+$ such that $b=cac'$. By definition $\le $ is always reflexive and transitive. 
For a Noetherian monoid, like $G^+$, $\le $ is a partial order, i.e. also antisymmetric. Indeed, assume there exist $a, b \in G^+$ such that $a \le b$ and
$b \le a$, i.e., there exist $c_1, c_2, c_1', c_2' \in G^+$ such that $ b=c_1ac_2$ and $a=c_1'bc_2'$. Then repeated insertion, i.e.,
$b=c_1ac_2=c_1(c_1'bc_2')c_2=c_1c_1'(c_1ac_2)c_2'c_2=\ldots $ contradicts Noetherianity except for the case $a=b$. 

\subsection{Normal forms in Garside systems} \label{3.2}
Since $\Delta $ is quasi-central, we have for all $a, b \in G$ and 
$k\in \mathbb{Z}$, $a\preceq \Delta ^k \preceq b$ iff $b \succeq \Delta ^k \succeq a$ iff $a \le \Delta ^k \le b$. 
If $\Delta ^r \le a\le \Delta ^s$, we denote $a\in [r,s]$. The maximal $r\in \mathbb{Z}$ and the minimal $s\in \mathbb{Z}$ such that $a\in [r,s]$ are called the {\it infimum} and the {\it supremum} of $a$, denoted by $\inf (a)$ and $\sup (a)$, respectively. \\
In Garside systems there exist natural normal forms. For every $a\in G^+$ there exists an unique $l\in \mathbb{N}$ and a unique decomposition $a=s_1\cdots s_l$
where $s_i=\Delta \wedge (s_i\cdots s_l)\in S$ and $s_l\ne 1$. This decomposition is called the {\it left normal form} of $a\in G^+$ in $(G,S)$, and the number $l$ is called the {\it canonical length} or {\it gap}, denoted $cl(a)$.
Using the decomposition $a=\Delta ^{\inf(a)}\overline{a}$ for some unique $\overline{a}\in G^+$ we get the {\it left $\Delta$-normal form} of $a$ in $(G,S)$.
Note that $\overline{a}_1=\overline{a}\wedge \Delta $ is a proper divisor of $\Delta $. We define $cl(a)=cl(\overline{a})$. \\
Furthermore, we can associate with every element $g\in G$ an unique pair $(a,b)\in G^+\times G^+$ such that $g=a^{-1}b$ and $a\wedge b=1$.
The left normal forms $a=t_1\cdots t_m$ and $b=s_1\cdots s_l$ give rise to the decomposition $g=(t_m^{-1}\cdots t_1^{-1})(s_1\cdots s_l)$ which we call the {\it left fractional normal form} of $g$ in $(G,S)$. {\it Right fractional normal forms} are defined symmetrically. \\
For the sake of brevity, in the sequel we often simply refer to a Garside group $G$ rather than a Garside system $(G, S)$ with $S={\rm Div}(\Delta )$.
 
\subsection{Garside structures for the braid groups} \label{GarsideBraids}

For any Garside system $(G, {\rm Div}(\Delta )$ and every $k\in \mathbb{N}$, also $(G, {\rm Div}(\Delta ^k)$ is a Garside system.
This leads us to the following definition.

\begin{defi} \label{MinGarsSys} \hspace{-0.15cm}{\bf .}
A pair $(G,S)$ is called a {\sc Minimal Garside system} if $(G, S)$ is a  Garside system and there exists no proper subset $S'$ of $S$ such that $(G, S')$
is also a Garside system. 
\end{defi}

$n$-strand braid groups admit at least two minimal Garside structures. The first being the classical Garside structure $(B_n, {\rm Div}(\Delta _n))$
where the partial orders $\preceq $ and $\succeq $ (defining the divisors of $\Delta $) are induced by the monoid generated by the atoms $\sigma _1, \ldots , \sigma _{n-1}$. Topologically, inside an Artin generator $\sigma _k$ ($1\le k < n$) the $(k+1)$-th strand crosses over the $k$-th strand, i.e., it is
an half twist on 2 strands. The Garside element $\Delta _n = \sigma _1 (\sigma _2\sigma _1)\cdots (\sigma _{n-1}\cdots \sigma _2\sigma _1)$  describes an half
twist on all $n$ strands. It has $n!$ positive divisors which are in one-to-one correspondence to permutations. Indeed, there exists a lattice isomorphism
from $( {\rm Div}(\Delta _n), \prec)$ to $S_n$ equipped with the weak order, where the fundamental braid $\Delta _n$ maps to the longest element. \\
The presentation with respect to the set of atoms is the classical Artin presentation \cite{\Ar}
\[ B_n=\left\langle \sigma _1,\ldots ,\sigma _{n-1} 
\left|
\begin{array}{ll} \sigma _i\sigma _j=\sigma _j\sigma _i & \forall |i-j|>1, \\ 
\sigma _i\sigma _{i+1}\sigma _i=\sigma _{i+1}\sigma _i\sigma _{i+1} & \forall i=1,\ldots ,n-2 \end{array} \right.\right\rangle.  \] 

The second, somehow dual (see \cite{\Be}), Garside structure is induced by the {\it dual monoid} $B_n^{+*}$ generated by the so-called {\it band generators}. Topologically, inside the band generator $a_{st}$ the $t$-th strand crosses over the $s$-th strand 'behind' the strands $s+1,\ldots ,t-1$, i.e.,
here the atoms correspond to all tranpositions (reflections) rather than only the nearest-neighbor transpositions (simple reflections).
W.r.t. this generating set $B_n$ admits the following complete set of relations \cite{\BKLb}.
\begin{eqnarray*}
a_{rs}a_{st}=a_{rt}a_{rs}=a_{st}a_{tr}, && 1\le r<s<t\le n, \\
a_{st}a_{qr}=a_{qr}a_{st}, && (t-r)(t-q)(s-r)(s-q)>0.
\end{eqnarray*}
We call this presentation the {\it dual presentation}. \\
The minimal Garside system induced by the dual braid monoid is $(B_n, {\rm Div}(\delta _n))$ with $\delta _n=a_{12}a_{23} \cdots a_{n-1,n}$.
Here the dual Garside element corresponds to a Coxeter element in the Weyl group. The simple elements are characterized by non-crossing partitions, i.e.,
we have $|{\rm Div}(\delta )|=\frac{1}{n+1} {2n\choose n}$, the $n$-th Catalan number.  \\
It is conjectured that for $n\ge 4$ $(B_n, \Delta _n)$  and $(B_n, \delta _n)$ are (up to automorphisms) the only minimal Garside systems for $B_n$. \\
The 3-strand braid group is exceptional as it admits further minimal Garside structures.

\begin{prop} \label{B3GarsSys} \hspace{-0.15cm}{\bf .}
(a) Recall that $B_3 \equiv T(2,3)=\langle a,b \mid a^2=b^3 \rangle $ where $a=\Delta _3$ and $b=\delta _3$ and consider the set $S_T=\{1,a,b,b^2,a^2\}$ of positive divisors of $a^2$ w.r.t the torus knot monoid $T^+=T^+(2,3)$. Then $(B_3, S_T)$ is a minimal Garside system.  \\
(b) $B_3$ also admits the non-standard presentation  $\langle \sigma ,b \mid \sigma b \sigma =b^2 \rangle $ where $\sigma =\sigma _1$ and $b=\delta _3$.
Consider the set $\{1,\sigma , b, \sigma b, b^2, b\sigma , b\sigma b, b^3\}$ of divisors of $b^3$ w.r.t. the monoid defined by that presentation.
Then $(B_3, {\rm Div}(b^3))$ is a minimal Garside system. 
\end{prop}

\subsection{Garside systems of spindle type} \label{3.4}

\begin{defi} \label{RigidGarsSys} \hspace{-0.15cm}{\bf .}
A Garside system $(G,S)$ with $S={\rm Div}(\Delta )$ is {\sc of spindle type} iff for every
pair of elements $a, b \in S$ we have that $a \wedge b$ and $a \, \tilde{\wedge } \, b$ lie in $\{ a, b, 1 \}$.
\end{defi}
In particular, if neither $a\preceq b$ nor $b\preceq a$ then $a\wedge b =1$, and if neither $a\succeq b$ nor $b\succeq a$ then $a \, \tilde{\wedge } \, b=1$.
\begin{coro} \label{JoinRigidGarsSys} \hspace{-0.15cm}{\bf .}
Let $(G,S)$ with $S={\rm Div}(\Delta )$ be a Garside system of spindle type. For every
pair of elements $a, b \in S$ we have that $a \vee b$ and $a \, \tilde{\vee } \, b$ lie in $\{ a, b, \Delta \}$.
\end{coro}

{\sc Proof.} The following formula holds for all $a, b,c, d\in G^+$ (see \cite{\Pb}):
\[ ac=bd=(a \vee b)(c \, \tilde{\wedge } \, d)=(a \wedge b)(c \, \tilde{\vee } \, d). \]
In particular, if $\Delta =a\partial (a)=b\partial (b)=\partial ^{-1}(c)c=\partial ^{-1}(d)d$ then 
\[ \Delta =(a \vee b)(\partial (a) \, \tilde{\wedge } \, \partial (b))=(\partial ^{-1}(c) \wedge \partial ^{-1}(d))(c \, \tilde{\vee } \, d). \]
We conclude that if the (left and right) gcd's are 1 then the lcm's equal $\Delta $.
And if $a \wedge b$ (or $a \, \tilde{\wedge } \, b$) lie in $\{ a, b \}$, then $a$ and $b$ are comparable w.r.t the partial order in question, and we
conclude that also $a \vee b$ (or $a \, \tilde{\vee } \, b$) lies in $\{ a, b\}$. \quad $\Box $  \\

This notion (with some additional height condition) was introduced by Picantin as {\it monoide de type fuseau} in \cite{\Pa}. \\
Note that all  Garside systems of spindle type are minimal. 

\begin{defi} \label{RigidGarsGP} \hspace{-0.15cm}{\bf .}
A group $G$ is a {\sc Garside group of spindle type} iff it admits a  Garside system  of spindle type $(G,{\rm Div}(\Delta ))$.
A monoid $M$ is {\sc  Garside monoid of spindle type} iff it is generated by $S={\rm Div}(\Delta )$ for some Garside system  of spindle type $(G,S)$.
\end{defi}

\noindent{\bf Examples.} (1) The 3-strand braid group $B_3$ is a  Garside group of spindle type. It admits the Garside monoids  of spindle type
$B_3^+$, $(B_3^*)^+$ and $T^+(2,3)$. However the Garside monoid $\langle \sigma ,b \mid \sigma b \sigma =b^2 \rangle ^+$ from  Proposition \ref{B3GarsSys} (b)
is not  of spindle type since $aba \wedge ba=b \notin \{aba, ba, 1\}$.  \\
(2) All {\it torus knot groups} $T(p,q)=\langle a,b \mid a^p=b^q \rangle $ are Garside groups of spindle type. \\
(2a) More general {\it torus-type groups} $T(p_1, \ldots ,p_k)=\langle a_1, \ldots , a_k \mid a_1^{p_1}=\ldots =a_k^{p_k} \rangle $ are 
Garside groups of spindle type. \\
(3) The Artin groups of type $I_2(m)$ (also called {\it dihedral Artin groups} $\mathcal{A}(I_2(m))$ are  Garside groups of spindle type. 
Indeed, two Garside structures are known - the classical one induced by the monoid $\mathcal{A}^+(I_2(m))\langle a,b \mid \underbrace{aba\cdots}_m=\underbrace{bab\cdots}_m \rangle $, and the dual induced by the monoid 
$\mathcal{A}^+(I_2^*(m))=\langle \{\sigma _i \}_{1\le i \le m} \mid \sigma _1\sigma _m=\sigma _m\sigma _{m-1} =\ldots =\sigma _2\sigma _1 \rangle $
where $\sigma _1=a$ and $\sigma _i=\underbrace{bab\cdots}_{i-1}=\underbrace{(bab\cdots)^{-1}}_{i-2}$ for $2\le i \le m$. \\
(4) The 3-strand pure braid group $P_3=\langle a ,b, c \mid  abc=bca=cab \rangle $ is a  Garside group of spindle type. Recall that for $n \ge 4$ the pure braid groups $P_n$
are not known to be Garside. In particular, the monoids induced by standard presentations are {\it not} Garside \cite{\Lee}. \\
(5)  $P_3$ is part of another infinite family of groups, namely the fundamental groups of the complement of a complex line arrangement whose graph is one multiple point of multiplicity $n \ge 2$. These fundamental groups admit presentations of the form
\[ \langle a_1,\ldots ,a_n \mid  a_na_{n-1}\cdots a_1=a_{n-1}\cdots a_1a_n=\ldots =a_1a_n\cdots a_2 \rangle , \]
and they are all Garside groups of spindle type with Garside element $\Delta =a_1 \vee \ldots \vee a_n=a_n \cdots a_1$.  
Indeed, these groups are the {\it pure dihedral Artin groups of rank n}, denoted by $\mathcal{PA}(I_2(n))$. The monoid $\mathcal{PA}^+(I_2(n))$ is a submonoid
of the dual dihedral monoid $\mathcal{A}^+(I_2^*(n))$ - the embedding is given by $a_i=\sigma _i^2$ for all $1\le i \le n$.  \\

\noindent{\bf Remark.} All Garside groups of spindle type given in the examples (except for example (2a)) are known to be {\it linear}, 
i.e. they admit a faithful linear representation.
This has been known for $B_3$ for many years. Indeed, for $n=3$ the Burau representation is faithful \cite{\MP}.
The dihedral Artin groups are of finite type. Faithful linear representations for finite type Artin groups are certain generalizations of the Lawrence-Krammer
representation \cite{\CW, \Di}.  \par
It can be shown that the torus knot groups are discrete subgroups of $\tilde{SL}_2(\mathbb{R})$ (see e.g. \cite{Tsa13}). \\

%\noindent{\bf Open problem.} Are the torus-type groups $T(p_1, \ldots ,p_k)$ linear?

\begin{prop} \label{PresRigidGarsSys} \hspace{-0.15cm}{\bf .}
Let $(G,S)$ with $S={\rm Div}(\Delta )$ be a Garside system of spindle type, and let $X=\{ a_1, \ldots , a_n \}$ be the set of atoms of $S$. 
Then $G$ admits the group presentation
\[ \langle a_1,\ldots ,a_n \mid  a_1 \underline{\partial (a_1)}=a_2 \underline{\partial (a_2)}=\ldots =a_n\underline{\partial (a_n)} \rangle.\] 
\end{prop}

{\sc Proof.}  By induction over $\preceq $, it is easy to show that every Garside group $G$ and every Garside monoid $G^+$ admits a {\it complemented presentation} in the sense of P. Dehornoy (see \cite{\DP}), i.e. $G=\langle X \mid a_i w_{ij}=a_j w_{ji} \,\, \forall i\ne j \rangle $ for some $w_{kl} \in X^*$.
with $\overline{w_{kl}}=a_k \backslash a_l$. In a Garside system of spindle type we have by Corollary \ref{JoinRigidGarsSys} $a_k \backslash a_l=\partial (a_k)$ and $a_k\vee a_l =\Delta $ for all $1\le k \ne l \le n$. \quad $\Box $  \\

\begin{prop} \label{PropsRigidGarsSys} \hspace{-0.15cm}{\bf .}
Let $(G,S)$ with $S={\rm Div}(\Delta )$ be a Garside system of spindle type and $X$ its set of atoms. \\
(1) If $\Delta \preceq b$ then a word $D\in X^*$ for $\Delta $ is an infix of any word $w\in X^*$ with $\overline{w}=b$, i.e. there exist
words $w_l, w_r\in X*$ s.t. $w\equiv w_lDw_r$. \\
(2) Let $b \in G^+ \setminus \Delta G^+$. Then there exists an unique word $w\in X^*$ representing $b$. \\   
(3) Let $\Delta ^q\preceq b$ for some $q\ge 0$. Then for every word $w\in X^*$ with $\overline{w}=b$ there exist $m, q_1,\ldots , q_m \in \mathbb{N}$,
$u_1, D_1 \ldots , u_m, D_m, u_{m+1} \in X^*$ such that $\sum _{i=1}^m q_i=q$, $\overline{D_i}=\Delta ^{q_i}$ $\forall i=1,\ldots m$, and
$w\equiv u_1 D_1 \cdots u_m D_m u_{m+1}$.
\end{prop}

{\sc Proof.}  Let $w, w'\in X^*$ s.t. $\overline{w}=\overline{w'}$. Then $w'$ can be transformed to $w$ by application of positive relations only, i.e. by
a sequence $\overline{w'}=\overline{w_1}= \ldots =\overline{w_m}=\overline{w}$ with $w_1, \ldots , w_m\in X^*$. Assume that a word for $\Delta  $ is not
an infix of $w$. Then, by Proposition \ref{PresRigidGarsSys}, no relation is applicable and we conclude that $w$ is the unique positive word representing $b$. Hence, we also get $\Delta \not \preceq \overline{w}$. Thus we have proven (1) and (2). \\
The proof of (3) is rather technical and relies on (1): Consider a decomposition $w\equiv u_1 D_1 \cdots u_m D_m u_{m+1}$
such that $\overline{D_i}=\Delta ^{q_i}$ $\forall i=1,\ldots m$, and $\sum _{i=1}^m q_i$ is maximal among all such decompositions of the word $w$.
Assume that $\sum _{i=j}^m q_i <q$. We will show that this leads to a contradiction.
Define words $u'_j \equiv \tau ^{\sum _{i=j}^m q_i} (u_j)$ for all $1 \le j \le m+1$, in particular $u'_{m+1} \equiv u_{m+1}$.
Since 
\[ \overline{w}=\overline{u_1 D_1 \cdots u_m D_m u_{m+1}} = \Delta ^{\sum _{i=j}^m q_i}\overline{u'_1} \cdots \overline{u'_{m+1}} \succeq \Delta ^q, \]
we conclude that $\overline{u'_1} \cdots \overline{u'_{m+1}} \succeq \Delta $. By (1) it follows that an infix of the word $u'_1 \cdots u'_{m+1}$ represents $\Delta $, i.e. there are two cases. Either there exists an $i_0 \in \mathbb{N}$ s.t. $u'_{i_0}$ contains an infix representing $\Delta $. 
Then so does $u_{i_0}$, and we conclude the decomposition $w \equiv u_1 D_1 \cdots u_{i_0}D_{i_0} \cdots u_m D_m u_{m+1}$ has not maximal $\sum _{i=1}^m q_i$.
In the second case, there are numbers $i_0 < j_0 \in \mathbb{N}$ and  words $u'^{(l)}_{i_0}, u'^{(r)}_{i_0}, u'^{(l)}_{j_0}, u'^{(r)}_{j_0}$ s.t.  
$u'^{(r)}_{i_0}$ is a postfix of $u'_{i_0}$, i.e. $u'_{i_0} \equiv u'^{(l)}_{i_0} u'^{(r)}_{i_0}$,
$u'^{(l)}_{j_0}$ is a prefix of $u'_{j_0}$, i.e. $u'_{j_0} \equiv u'^{(l)}_{j_0}u'^{(r)}_{j_0}$,
and $\overline{u'^{(r)}_{i_0} u'_{i_0+1} \cdots u'_{j_0-1} u'^{(l)}_{j_0}} =\Delta $.
Decompose $u_{i_0} \equiv u^{(l)}_{i_0} u^{(r)}_{i_0}$ with $u^{(l),(r)}_{i_0} \equiv \tau ^{- \sum _{i=i_0}^{m} q_i} (u'^{(l),(r)}_{i_0})$
and $u_{j_0} \equiv u^{(l)}_{j_0} u^{(r)}_{j_0}$ with $u^{(l),(r)}_{j_0} \equiv \tau ^{- \sum _{i=j_0}^{m} q_i} (u'^{(l),(r)}_{j_0})$.
Then $u^{(r)}_{i_0} D_{i_0} u_{i_0 +1} D_{i_0+1} \cdots u_{j_0 -1} D_{j_0-1} u^{(l)}_{j_0}$ represents
\[ \Delta ^{\sum _{i=i_0}^{j_0-1} q_i} \tau ^{- \sum _{i=j_0}^{m} q_i} (\overline{u'^{(r)}_{i_0} u'_{i_0+1} \cdots u'_{j_0-1} u'^{(l)}_{j_0}}) 
= \Delta ^{(\sum _{i=i_0}^{j_0-1} q_i) +1}. \]
Set $\tilde{D}_i \equiv D_i$ for $1 \le i \le i_0-1$, 
$\tilde{D}_{i_0} \equiv u^{(r)}_{i_0} D_{i_0} u_{i_0 +1} D_{i_0+1} \cdots u_{j_0 -1} D_{j_0-1} u^{(l)}_{j_0}$,
and $\tilde{D}_j \equiv D_{j+j_0-i_0}$ for $i_0+1 \le j \le m-j_0+i_0$. Then, for all $1\le i \le m-j_0+i_0$, 
$\overline{\tilde{D_i}}=\Delta ^{\tilde{q}_i}$ for some $\tilde{q}_i \in \mathbb{N}$. Indeed, $\tilde{q}_i=q_i$ for $1\le i \le i_0-1$,
$\tilde{q}_{i_0}=(\sum _{i=i_0}^{j_0-1} q_i) +1$, and $\tilde{q}_j=q_{j+j_0-i_0}$ for $i_0+1\le j \le m-j_0+i_0$.
In addition, set $\tilde{u}_i=u_i$ for $1\le i \le i_0-1$,
$\tilde{u}_{i_0}=u^{(l)}_{i_0}$, $\tilde{u}_{i_0+1}=u^{(r)}_{j_0}$, and $\tilde{u}_j=u_{j+j_0-i_0-1}$ for $i_0+2\le j \le m-j_0+i_0+1$,
and define $\tilde{m}:=m-j_0+i_0$.
Then by construction 
\[ \tilde{u}_1 \tilde{D}_1 \cdots \tilde{u}_{\tilde{m}} \tilde{D}_{\tilde{m}} \tilde{u}_{\tilde{m}+1} \equiv  u_1 D_1 \cdots u_m D_m u_{m+1} \equiv w \]
with $\sum _{i=1}^{\tilde{m}} \tilde{q}_i= (\sum _{i=1}^{m} q_i)+1$ in contradiction with the maximality of $\sum _{i=1}^{m} q_i$.
\quad $\Box $  \\

\subsection{Logspace normal forms for Garside systems of spindle type}  \label{3.5}
\subsubsection{Garside normal form}  \label{3.5.1}

\begin{prop} \label{LNFRigidGarsSys} \hspace{-0.15cm}{\bf .}
Let $(G,S)$ with $S={\rm Div}(\Delta )$ be a Garside system of spindle type. If $G$ is a linear group, then the (left) Garside normal form 
of any element $b\in G$ can be computed by a logspace transducer.  
\end{prop}

{\sc Proof.} 
Let $X=\{a_1, \ldots , a_n\}$ be the set of atoms of the given Garside system of spindle type $(G,{\rm Div}(\Delta ))$.
Assume that our input element is given as a word $w\in (X^{\pm })^*$.  Recall that transforming a word $v$ over another alphabet $Y$ to a word $w\in (X^{\pm })^*$ can be done in logarithmic space. We will construct explicitly an $L$-computable function $f$ over $(X^{\pm })^*$ that (on input $w$) returns its
(left) Garside normal form, i.e., a tuple $(p, s_1, \ldots , s_l) \in \mathbb{Z} \times (S\setminus \{1, \Delta \})^*$ such that $\overline{w}\stackrel{LNF}{=}\Delta ^p s_1\cdots s_l$. 
Indeed, $f$ is the composition of four $L$-computable functions $f_a, f_b, f_c, f_d$. \\
(a) The function $f_a: (X^{\pm })^* \longrightarrow (X^{\pm })^*$ returns a freely reduced word $w_a$ s.t $\overline{w_a}=\overline{w}$. Recall that free reduction in the rank $n$ free group $\langle a_1, \ldots ,a_n$ can be done in Logspace. Actually, step (a) is not really necessary, but we include it for convenience. \\
(b)  The function $f_b: (X^{\pm })^* \longrightarrow \mathbb{N} \times X^*$ maps $w_a$ to $(k, w_b)$ s.t. $\overline{w_a}=\Delta ^{-k}\overline{w_b}$. Here $k\ge 0$ is the number of occurences of negative generators in the word $w_a$. Let $w_a \equiv u_1a_{i_1}^{-1}u_2a_{i_2}^{-1}\cdots u_ka_{i_k}^{-1}u_{k+1}$ with
$u_j \in X^*$ for $1\le j \le k+1$. Replacing, for $1\le j \le k$, each $a_{i_j}^{-1}$ by $\partial (a_{i_j}) \Delta ^{-1}$, and sliding the $\Delta ^{-1}$'s to
the left, leads to the following word (equivalent to $w_a$):
\[ \Delta ^{-k} \tau ^{-k}(u_1\partial(a_{i_1})) \tau ^{-k+1}(u_1\partial(a_{i_2})) \cdots \tau ^{-1}(u_1\partial(a_{i_k})) u_{k+1} \equiv \Delta ^{-k}w_b. \]  
It is easy to compute $f_b(w_a)$ in logarithmic space. First scan the word $w_a$ and count the number of negative occurences of atoms, and output (and store) that number $k$ whose bitlength is $O(\log |w_a|)$. Set a counter $j:=k$. Scan again the word $w_a \equiv x_1\cdots x_i \cdots x_{|w_a|}$ and do the following:
\begin{itemize}
 \item[-] If $x_i$ is an atom, then output $\tau ^{-j}(x_i)$. 
 \item[-] If $x_i$ is the inverse of an atom, then output the letters of the unique word representing $\tau ^{-k}(\partial(x_i^{-1})) 
                 \in S\setminus \{1, \Delta \}$, and decrement $j:=j-1$. 
\end{itemize}
(c) The function $f_c: \mathbb{N} \times X^* \longrightarrow \mathbb{N} \times X^* \times \mathbb{N}$ maps $(k,w_b)$ to $(k,w_c, q)$ 
such that $\overline{w_b}= \overline{w_c}\Delta ^q$. 
and $q\ge 0$ is maximal such that $\Delta ^q \preceq \overline{w_b}$. We conclude that $\inf (\overline{w})=q-k$.  
To find $q$ we have to extract the whole ``$\Delta $-content'' from the positive word $w_b$. \\
Since ``nesting'' might occur as e.g. in the word $babbaababb$ representing the 4-th power of the Garside element in $T(2,3)$, we cannot use Proposition \ref{PropsRigidGarsSys} (a) by removing one $\Delta $ after the other for a logspace algorithm.  
However we can utilize Proposition \ref{PropsRigidGarsSys} (c) which states
that there exists a decomposition $w\equiv u_1 D_1 \cdots u_m D_m u_{m+1}$ with 
$u_1, D_1 \ldots , u_m, D_m, u_{m+1} \in X^*$, $\overline{D_i}=\Delta ^{q_i}$ $\forall i=1,\ldots m$, and $\sum _{i=1}^m q_i=q$.
Indeed, the following algorithm utilizes such a decomposition with $m$ minimal. \\ 
For a monoid element $b \in G^+$, we define the so-called ``norm'' by $||b ||=\max \{|w| \mid w\in X^*, \, \overline{w}=b \}$,
and for any element $b\in G$, denote the geodesic length (w.r.t. the generating set $X$) $\ell (b)=\ell _X(b)=
\min \{|w| \mid w\in (X^{\pm })^*, \, \overline{w}=b \}$. Let $w_b \equiv x_1\cdots x_i \cdots x_{|w_b|}$ with $x_i\in X$.

\begin{algorithm}
\caption{Function $f_c: \mathbb{N} \times X^* \longrightarrow  \mathbb{N} \times X^* \times \mathbb{N}$.}
%\SetLine
\KwIn{$(k, w_b) \in \mathbb{N} \times X^*$.}
\KwOut{$(k, w_c, q) \in \mathbb{N} \times X^* \times \mathbb{N}$}
Print $k$ on output tape\;
Initialize counter $q:=0$\;
Initialize $i:=1$\;
\While{$i \le |w_b|$}{
 \For{$j:=|w_b|$ {\bf to} 1 {\bf by} -1}{
   $q_{\min }:=\lceil (j-i+1)/ ||\Delta || \rceil $;
   $q_{\max }:=\lfloor (j-i+1)/ \ell (\Delta ) \rfloor $\;
   \For{$q_1:=q_{\min }$ {\bf to} $q_{\max }$}{
      Check whether $x_i\cdots x_j \Delta ^{-q_1}\stackrel{!}{=} 1$ calling the WP-oracle \;
      \If{$x_i\cdots x_j \Delta ^{-q_1} = 1$}{$i:=j$; $q:=q+q_1$\; {\bf break} $j$\;}
   }   
 }
 \lIf{$x_i\cdots x_j \Delta ^{-q_1} \ne 1$ {\rm for all} $q_1\in [q_{\min }(i,j), q_{\max }(i,j)]$ $\forall j: i\le j \le k$}
      {Print $\tau ^{-q}(x_i)$ on output tape\;}
 $i:=i+1$\;
}
Print $q$ on output tape.
\label{fc}
\end{algorithm}

Algorithm \ref{fc} checks for each letter $x_i$ whether it is the first letter in some word representing a $\Delta $-power.
Only letters $x_i$ not belonging to some $\Delta $-power word (according to that procedure) are printed, but as $\tau ^{-q}(x_i)$ 
(here we mean the current $q$-value during the runtime of the algorithm) because we have to move $\Delta ^q$ through such letters to the right.
For checking whether $x_i\cdots x_j \Delta ^{-q_1}\stackrel{?}{=} 1$ we need to call the oracle for the word problem (WP)
which is so far known to be logspace computable only for linear (!) Garside groups of spindle type. \\
For all $D\in X^*$ with $\overline{D}=\Delta $, we have $\ell (\Delta ) \le |D| \le ||\Delta || $ which implies
$q_1 \ell (\Delta ) \le |D^{(q_1)}| \le q_1 ||\Delta || $ for all words  $D^{(q_1)}\in X^*$ representing $\Delta ^{q_1}$.
We conclude that if $x_i\cdots x_j =\Delta ^{q_1}$ then \[ (j-i+1)/||\Delta || \le q_1 \le  (j-i+1)/ \ell (\Delta ), \]
and we may choose integer test values for $q_1$ from that interval $[q_{\min }, q_{\max }]$.  \\
(d) The last function $f_d: \mathbb{N} \times X^* \times \mathbb{N} \longrightarrow \mathbb{Z} \times (S\setminus \{ 1, \Delta \})^*$ maps $(k, w_c, q)$ to $(q-k, s_1, \ldots , s_l)$ such that $\Delta ^{-k} \overline{w_c} \Delta ^q \stackrel{LNF}{=} \Delta ^{q-k} s_1\cdots s_l$. 
Since $\inf (\overline{w_c})=q-k$ $\overline{w_c}$ lies in $G^+ \setminus \Delta G^+$ and $w_c$ is according to \ref{PropsRigidGarsSys} the unique
representative of $\overline{w_c}$ and it is (in some sense) its left and right Garside normal form (which coincide).
The same holds for its automorphic image $\tau ^q(w_c)$ - recall that $\Delta ^{-k} \overline{w_c}\Delta ^q=\Delta ^{q-k} \tau ^q(\overline{w_c})$. \\
Therefore, we only have to read off the simple elements from the word $\tau ^q(w_c)$. Let $w_c=x_1 \cdots x_{|w_c|}$ with $x_i \in X$ for all $i$.
Furthermore, let $\underline{S}$ be the set of unique representatives of the elements of $S\setminus \{ 1, \Delta \}$.
We scan the word $w_c$ for consecutive subwords, and we print $s\in S$ (or equivalently its unique representative $\underline{s}\in \underline{S}$) 
on the output tape iff  $\tau ^q(x_i)\cdots \tau ^q(x_j)=s$ ($i<j<|w_c|$), but the word $x_i \cdots x_jx_{j+1}$ lies not in $\underline{S}$. Of course, we also print the last simple element $s_l=\tau ^q(x_{i_l}) \cdots \tau ^q(x_{|w_c|})$ (for some minimal $i_l$). 
\quad $\Box $  \\

\subsubsection{Geodesic normal form} \label{3.5.2}

\begin{defi} \label{homogPres} \hspace{-0.15cm}{\bf .}
Let $(G,S)$ be a Garside system and $X$ the set of atoms of $S$. $(G, S)$ is called a Garside system with {\sc homogenous presentation} iff $G$ (or/and $G^+$) 
admits a presentation of the form $\langle X \mid R^{(l)}_1=R^{(r)}_1, \ldots R^{(l)}_m=R^{(r)}_m \rangle $ with $R^{(l)}_i, R^{(r)}_i \in X^*$
and $|R^{(l)}_i|=|R^{(r)}_i|$ for all $1 \le i \le m$.
\end{defi}

For a Garside system $(G,S)$ with homogeneous presentation all words $w\in X^*$ representing a monoid element $p\in G^+$ have the same word length
$|w|=\ell (p)=||p||$.

\begin{defi} \label{RedOp} \hspace{-0.15cm}{\bf .}
Let $(G,S)$ be a Garside system  with homogeneous presentation and $X$ the set of atoms of $S$. We define a Reduction operations 
${\rm red}, {\rm Red}: (X^{\pm })^* \longrightarrow (X^{\pm })^*$. Consider the unique decomposition of a word $w\in (X^{\pm })^*$ into a form
$w\equiv \Delta ^r w_1 \cdots w_s$ with $w_i \in (S\setminus \{1, \Delta \})^{\pm }$ such that $s$ is minimal. Then we define ${\rm red}(w)$ as follows: \\   
(1) If $r\ge 0$ or $w_i \preceq 1$ then set ${\rm red}(w) \equiv w$. \\
(2) Otherwise, choose a $w_k \succeq 1$  whose wordlength is maximal among $\{ w_i \in (0,1) \mid 1\le i \le s$, and define
\[ {\rm red}(w) \equiv \Delta ^{r+1} \tau (w_1) \cdots \tau (w_{k-1}) w'_k w_{k+1} \cdots w_s \]
where $w'_k \in (X^-)^*$ is a word of minimal length s.t. $\Delta \overline{w'_k}=\overline{w_k}$.  \\
Now, we set ${\rm Red} (w) \equiv {\rm red}^{|r|}(w)$.
\end{defi}

Note that although neither ${\rm red}(w)$ nor ${\rm Red}(w)$ is uniquely determined from a
given decomposition of $w$, the word-length $|{\rm Red}(w)|$ is uniquely determined and $|{\rm Red}(w)| \le |w|$. \\

\noindent{\bf Remark.} We restricted Definition \ref{RedOp} to Garside systems  with homogeneous presentation, because for
Garside systems  with non-homogeneous presentation Reduction might increase word length, i.e, $|{\rm Red}(w)| \le |w|$ does not hold. 
Consider for example a torus knot group $\langle a,b \mid a^p=b^q  \rangle $ with $p+2 < q$. Here, reduction of the word $w\equiv a^{-p}b$ 
leads to the word ${\rm Red}(w)\equiv b^{-q+1}$ with $|{\rm Red}(w)|=q-1 > |w|=p+1$.

\begin{defi} \label{ShWP} \hspace{-0.15cm}{\bf .}
The {\sc shortest (or geodesic) word problem} (w.r.t. the generating set $X$) is
the following problem. Given a word $w\in (X^{\pm })^*$, find $\ell (w)=\ell _X(w)$ and return a geodesic representative, i.e. a
word $w_g$ s.t. $|w_g|=\ell(w)$ and $\overline{w_g}=\overline{w}$. 
\end{defi}

\begin{lem} \label{LemShWP} \hspace{-0.15cm}{\bf .}
Let $(G,S)$ be a Garside system of spindle type with homogeneous presentation and $X$ the set of atoms of $S$. Let $w_{LNF}$ be a word in LNF, i.e. $w_{LNF} \equiv D^p w_1\cdots w_l$ with $D, w_1, \ldots , w_l \in X^*$ such that $\overline{D}=\Delta $, $\overline{w_i}=s_i \in S$ 
and $\Delta ^ps_1 \cdots s_l$ is a (left) Garside normal form. \\
Then, for all words $w\in (X^{\pm })^*$ such that $\overline{w}=\overline{w_{LNF}}$, 
\[ |{\rm Red}(w_{LNF})| \le |w|. \]  
\end{lem}

{\sc Proof.}  We start from a decomposition of $w$, say $w \equiv x_1 \cdots x_{|w|}$,
where $x_i$ is either an atom or an inverse of an atom. Recall that $|{\rm Red}(w)| \le |w|$.
If there is a factor, say $x_k$, which is an inverse of a generator, change $w$ to
\[ w' \equiv \Delta ^{-1} \tau ^{-1}(x_1) \cdots \tau ^{-1}(x_{k-1}) w'_k x_{k+1} \cdots x_{|w|} \]
where $w'_k$ is a word of minimal length such that $\Delta ^{-1}\overline{w'_k}=x_k$, and therefore $w'_k \in X^*$. Then clearly,
$|{\rm Red}(w')| \le |{\rm Red}(w)|$.
Continue this process until we get a decomposition $\Delta ^r c_1 \cdots c_m$ with $c_i \in S\setminus \{1, \Delta \}$.
To prove the lemma, it suffices to show that the value $|{\rm Red}(\cdot )|$ does not increase
for all decompositions that arise during the transformation of $\Delta ^r c_1 \cdots c_m$ into the left-canonical
decomposition  $\Delta ^p s_1 \cdots s_l$. Let $c_ic_{i+1}$ be transformed to  $c'_ic'_{i+1}$. It suffices to show that
\[ \max \{ |c'_i|, |c'_{i+1}| \} \ge \max \{ |c_i|, |c_{i+1}| \} \quad (*) \]
which implies $|{\rm Red}( \Delta ^qc'_ic'_{i+1})| \le |{\rm Red}(\Delta ^q c_ic_{i+1})|$ for all $q \in \mathbb{Z}$.
Note that in Lemma 5.1. of \cite{\KKL} $(*)$ is proven for all homogeneous Garside systems with $|| \Delta ||=3$ (which implies
$1\le |c_i|, |c_{i+1}| \le 2$).  \\
By definition of leftgreedy decomposition, we have $c_ic_{i+1}=c'_ic'_{i+1}$ with $c'_{i+1}=\Delta \wedge (c_ic_{i+1})$, i.e. $\forall e: e \preceq c_ic_{i+1}
\Rightarrow e\preceq c'_i$, in particular $c_i \preceq c'_i$. Write $c'_i=c_id$. Since $G$ is  of spindle type, this implies $d \preceq \partial (c_i)$.
By left cancellativity, we conclude that $c_{i+1}=dc'_{i+1}$. In case of $d=1$ we get $c'_i=c_i$ and $c'_{i+1}=c_{i+1}$, and $(*)$ holds trivially.
Now, assume $1 \ne d \preceq \partial (c_i)$. Recall that $c_{i+1}$ is simple, i.e. $c_{i+1}=dc'_{i+1} \preceq \partial (c_i) \partial (\partial (c_i))=\Delta $.
Here we conclude that either $d= \partial (c_i)$ and $c'_{i+1} \preceq \partial (\partial (c_i))$, or $d \ne \partial (c_i)$ and $c'_{i+1}=1$.
In both cases we have $|c'_i| \ge \max \{ |c_i|, |c_{i+1}| \}$, validating $(*)$.
\quad $\Box $  \\

\begin{theo} \label{ShWPrigidGS} \hspace{-0.15cm}{\bf .}
Let $(G,S)$ be a  Garside system of spindle type with homogeneous presentation and $X$ the set of atoms of $S$. 
If $G$ is a linear group then we can solve the geodesic word problem w.r.t $X$ with a logspace transducer. Indeed, we may always find some unique $w_g$, i.e.
a {\sc geodesic normal form} in Logspace. 
\end{theo}

{\sc Proof.} Given an instance word $w \in (X^{\pm })^*$, recall that Lemma \ref{LemShWP} already solves the shortest WP w.r.t. to $\ell _X(\cdot )$ by providing the geodesic $w_g={\rm Red}(LNF(w))$. It remains to make the reduction procedure unique and to show that it can be accomplished in logspace. First, 
according to Proposition \ref{LNFRigidGarsSys}, computing the (unique) left normal form $\Delta ^p s_1\cdots s_l$ can be done with logarithmic space.
Let $w_{LNF}=D^p w_1 \cdots w_l$ be the (unique for Garside systems {\it of spindle type}) LNF word, i.e., 
$D, w_1, \ldots , w_l \in X^*$ s.t. $\overline{D}=\Delta $ and $\overline{w_i}=s_i$ for all $1\le i \le l$.
Now, if $p\ge 0$, then clearly ${\rm Red}(w_{LNF})=w_{LNF}$. 
Also if $p<0$ and $|p|\ge l$, then ${\rm Red}(w_{LNF})=D^{p+l}\tau ^l(w'_1) \tau ^{l-1}(w'_2) \cdots \tau (w'_{l-1})w'_l$ is already unique,
where $w'_k \in (X^-)^*$ is a word of minimal length s.t. $\Delta \overline{w'_k}=\overline{w_k}$ (see Def. \ref{RedOp}). \\
Therefore, one may assume that $p<0$ and $|p|<l$. Define a threshold value as the greatest number $t\in \mathbb{N}$ s.t. at least $p$ of the words $w_k$ have length $\ge t$, i.e $|\{ w_k \mid 1\le k\le l, \,\, |w_k|\ge t \}|\ge p$. This value $t$ can be computed in logspace by running though all values 
$t=1,\ldots , ||\Delta ||-1$, and for each $t$ checking whether $|\{ w_k \mid 1\le k\le l, \,\, |w_k|\ge t \}|\ge p$. Clearly this can be done with the help
of two logspace counters. Note that one may improve that procedure by running only through all possible values of 
$||s||=\ell (s)$ for all $s\in S\setminus \{1, \Delta \}$. \\
Now, we make our Red-operation unique my applying the red-operation to the $|p|$ leftmost words $w_k$ with $|w_k|\ge t$, i.e let $1\le i_1 < \ldots < i_p$
with $|w_{i_k}|\ge t$ $\forall 1\le k\le p$ and $\sum _{k=1}^p i_k$ minimal. Then one may decompose (in logspace) the LNF word as
$w_{LNF} \equiv D^p u_0 w_{i_1} u_1 \cdots w_{i_p}u_p$ for some $u_0, \ldots , u_p\in X^*$. 
Recall that $w'_k \in (X^-)^*$ is a word of minimal length s.t. $\Delta \overline{w'_k}=\overline{w_k}$.
Then we output the following geodesic word 
\[  w_g\equiv \tau ^{|p|}(u_0) \tau ^{|p|-1}(w'_{i_1} \cdot u_1) \cdots \tau (w'_{i_{|p|-1}} \cdot u_{|p|-1}) \cdot w'_{i_{|p|}}\cdot u_{|p|}, \]
which can easily be accomplished with logarithmic space only.
\quad $\Box $  \\

\section{Conjugacy in Garside groups of spindle type} \label{IV}

\subsection{Cycling and Decycling}

Let us recall some further definitions and facts in Garside theory.
\begin{defi} \label{Order} \hspace{-0.15cm}{\bf .}
Let $(G,S)$ with $S={\rm Div}(\Delta )$ be a Garside system. Recall that $\tau : G \rightarrow G$ is defined by $a \mapsto \Delta ^{-1}a\Delta $.
The {\sc order of a Garside system}, denoted by ${\rm ord} (G,S)$, is the minimal number $n\in \mathbb{N}$
such that $\tau ^n={\rm id}$. Analogeously, one may define for any $b\in G$ the {\sc order of the element} ${\rm ord}  (b)$ (w.r.t. that Garside system)
as the  minimal  number $n\in \mathbb{N}$ such that $\tau ^r(b)=b$. By definition, ${\rm ord}  (b)$ is always a divisor of ${\rm ord}  (G,S)$.
\end{defi}

In particular, the following definitions are fundamental for the study of conjugacy in Garside groups.
\begin{defi} \label{SSS} \hspace{-0.15cm}{\bf .}
Let $(G,S)$ with $S={\rm Div}(\Delta )$ be a Garside system and $b$ an element in $G$.  We define the {\sc summit infimum}, the {\sc summit supremum} and the {\sc summit canonical length} (or {\sc summit gap}) as the maximal (or minimal) possible value of these quantities inside the conjugacy class of $b$, i.e.,
\[ {\inf }_s(b):=\max \{ \inf (b') \mid b' \sim b \}, \quad {\sup }_s(b):=\min \{ \sup (b') \mid b' \sim b \},\]
and $cl_s(b):={\sup }_s(b)-{\inf }_s(b)$.
For $p\le q\in \mathbb{N}$, define intervals $[p,q]:=\{ a\in G \mid \Delta ^p \preceq a \preceq \Delta ^q\}$. 
Then the {\sc Super Summit Set} $SSS(b)$ of an element $b \in G$ is defined as the intersection of the conjugacy class $C(b)$ with the interval $[{\inf }_s(b), {\sup }_s(b)]$.
\end{defi}

Summit infimum, supremum and gap are conjugacy invariants. The SSS is a full invariant, i.e., $b_1 \sim b_2$ if and only if $SSS(b_1)=SSS(b_2)$.
A representative $\tilde{b}$ in $SSS(b)$ may be found by iterative application of a finite number of {\it cycling} and {\it decycling} operations $c, d$.
Let $\Delta ^ps_1 \cdots s_l$ be the LNF of $b$. Then
\[ c(b):=\Delta ^p s_2\cdots s_l \cdot \tau ^{-p}(s_1), \, \, {\rm and} \quad  d(b):=\Delta ^p \tau ^p(s_l)\cdot s_1 \cdots s_{l-1}. \]

For a general Garside system $(G, {\rm Div} (\Delta ))$, according to an improvement \cite{\BKLa} of the {\it cycling theorem} \cite{\EM} one must ``cycle'' (resp. ``decycle'') at most $||\Delta ||-1$ times in order to either increase $\inf $ (resp. decrease $\sup $) or to be sure that it is already maximal (resp.
minimal) for the given conjugacy class. For Garside systems of spindle type the situation is even much simpler.
Elements in the SSS can be characterized as follows.

\begin{prop} \label{SSSrigid} \hspace{-0.15cm}{\bf .}
Let $(G,S)$ with $S={\rm Div}(\Delta )$ be a Garside system of spindle type and $b$ an element in $G$ with left Garside normal form $\Delta ^p s_1 \cdots s_l$.  
Then \[  b \in SSS(b) \quad if \,\, and \,\, only \,\, if \quad \tau ^{-p}(s_1) \wedge \partial (s_l)=1. \]
In this case (if $\tau ^{-p}(s_1) \wedge \partial (s_l)=1$) we call the element $b$ {\sc rigid}. \\
In particular, if $b$ is a not rigid then
\[ (a) \quad \inf (d(b))=\inf (b) +1 \quad {\it or} \quad \sup (d(b))=\sup (b) -1. \]
\end{prop}
{\sc Proof.} ($\Rightarrow $): Let $s_1':= s_1 \wedge \tau ^p(\partial (s_l))$ and write $s_1=s_1's_1''$. 
The negation of the r.h.s is equivalent to $s_1'\ne 1$.  Since $G$ is of spindle type, we conclude that $s_1'=s_1$ or $s_1'=\partial (\tau ^p(s_l))$.
If $s_1'=s_1$ then $s_1\preceq \partial (\tau ^p(s_l))$, and therefore $s:=\tau ^p(s_l)s_1 \preceq \Delta $. We conclude that $d(b)=b^{s_l}=
\Delta ^p s s_2\cdots s_{l-1}$ has supremum $\sup (d(b))=p+l-1=\sup (b)-1$, and therefore $b \notin SSS(b)$. \\
Now, if $s_1'=\partial (\tau ^p(s_l))$ then $d(b)=\Delta ^{p+1} s_1'' s_2\cdots s_{l-1}$, i.e. $\inf (d(b))=\inf (b) +1$. Again we conclude that
$b \notin SSS(b)$. Obeserve that we proved here also the second claim. \\ 
($\Leftarrow $): Now let $b$ be rigid, i.e. $s_1 \wedge \tau ^p(\partial (s_l))=1$, and assume that $b\notin SSS(b)$. Then, according to the cycling theorem \cite {\EM}, after applying finally many cyclings/decyclings one should increase/decrease the infimum/supremum of $b$. Since $s_1 \wedge \tau ^p(\partial (s_l))=1 \Leftrightarrow \tau ^p(s_l)s_1 \wedge \Delta =\tau ^p(s_l)$ the pair $(\tau ^p(s_l),s_1)$ is left-greedy. Therefore $\Delta ^p \tau ^p(s_l)\cdot s_1 \cdots s_{l-1}$ is already the left (and right) Garside normal form of $d(b)$.

%That's the reason why we call $b$ rigid. Anyway, $d(b)$ has the same $\inf $ and $\sup $ as $b$.

 By induction over $k$ we may show that  
\[ d^k(b)=\Delta ^p \tau ^p(s_{l-k+1} \cdots s_l) s_1\cdots s_{l-k} \quad {\rm for} \quad 1\le k \le l.  \]
In particular, $d^l(b)=\Delta ^p \tau ^p(s_1 \cdots s_l)=\tau ^p(b)$. 
Let $r={\rm ord}(b) \mid {\rm ord} (G)$, then $d^{lcm (p,r)l/p}(b)=\tau ^{lcm (p,r)}(b)=b$ for $p>0$ and $d^l(b)=b$ for $p=0$. Since, the supremum can not increase by applying decyclings, we conclude
that $\sup d^k(b)=\sup (b)$ for all $k\in \mathbb{N}$.
Analogeously, one may show for cyclings that 
\[ c^k(b)=\Delta ^p  s_{k+1}\cdots s_l \tau ^{-p}(s_1 \cdots s_k) \quad {\rm for} \quad 1\le k \le l, \quad c^l(b)=\tau ^{-p}(b),  \]
and $c^{lcm (p,r)l/p}(b)=\tau ^{-lcm (p,r}(b)=b$.  
Since, the infimum can not decrease by applying cyclings, we conclude
that $\inf c^k(b)=\inf (b)$ for all $k\in \mathbb{N}$. Applying the cycling theorem leads to the conclusion that $b$ has already maximal infimum and 
minimal supremum inside its conjugacy class, i.e. $b\in SSS(b)$.
\quad $\Box $ \\

For a Garside groups of spindle type, if an element lies not inside its SSS, then we may find an element with bigger/smaller infimum/supremum by only one decycling.
A similar statement holds for cyclings. 

\begin{prop} \label{Decyc} \hspace{-0.15cm}{\bf .}
Let $(G,S)$ with $S={\rm Div}(\Delta )$ be a Garside system of spindle type. Consider a non-rigid element $b\stackrel{LNF}{=}\Delta ^p s_1 \cdots s_l$.  
If $\tau ^p(s_l)s_1 \ne \Delta $, then $d(b)\in SSS(b)$.
\end{prop}
{\sc Proof.} 
Again denote $s_1':= s_1 \wedge \tau ^p(\partial (s_l))$. Since $(G,S)$ is of spindle type and $s_1' \ne 1$ ($b$ is non-rigid), we have two cases to consider. \\ 
(a) Either we have $s_1'=s_1$ which implies $s_1 \preceq \partial (\tau ^p(s_l))$.
Since $\tau ^p(s_l)s_1 \ne \Delta $ we conclude that  $\tau ^p(s_l)s_1 \prec \Delta $. 
The left-greedy condition for the pair $(s_{l-1},s_l)$ is equivalent to $\tau ^p(s_l) \wedge \partial (\tau ^p(s_{l-1}))=1$.
We conclude that $\tau ^p(s_l)s_1 \wedge \partial (\tau ^p(s_{l-1}))=1$. Since the LNF of $d(b)$ is $\Delta ^p ss_2\cdots s_{l-1}$ with $s=\tau ^p(s_l)s_1$,
this is exactly the rigidity condition for $d(b)$. By Proposition \ref{SSSrigid} we conclude that $d(b)\in  SSS(b)$. \\
(b) The other case is $s_1'=\partial (\tau ^p(s_l))$. Hence $\Delta \preceq \tau ^p(s_l)s_1=\Delta s_1''$ with $s_1=s_1's_1''$, and
\[ 1 \ne s_1'' \preceq \partial (s_1')= \partial (\partial ( \tau ^p (s_l)))=\tau ^{p+1}(s_l). \]  
 The left-greedy condition for the pair $(s_{l-1},s_l)$ is equivalent to $\tau ^{p+1}(s_l) \wedge \partial (\tau ^{p+1}(s_{l-1}))=1$.
We conclude that $d \wedge \partial (\tau ^{p+1}(s_{l-1}))=1$ for all positive left divisors $d\preceq \partial (\tau ^{p+1}(s_{l-1}))$. 
In particular, this holds for $d=s_1''$, i.e.  $s_1'' \wedge \partial (\tau ^{p+1}(s_{l-1}))=1$
Since the LNF of $d(b)$ is $\Delta ^{p+1} s_1''s_2\cdots s_{l-1}$, this is exactly the rigidity condition for $d(b)$. 
By Proposition \ref{SSSrigid} we conclude that $d(b)\in  SSS(b)$.
 \quad $\Box $ \\

\begin{lem} \label{DecycPow} \hspace{-0.15cm}{\bf .}
Let $(G,S)$ with $S={\rm Div}(\Delta )$ be a  Garside system of spindle type. Consider a non-rigid element $b\stackrel{LNF}{=}\Delta ^p s_1 \cdots s_l$.  
Assume that the summit infimum is $\inf _s(b)=p+p_0$. Then, for $1\le k \le p_0$,
\[ d^k(b)=\left\{ \begin{array}{ll} \Delta ^{p+k} s_{k+1} \cdots s_{l-k}, & k<p_0 \\
\Delta ^{p+p_0} s_{p_0}'' s_{p_0+1}\cdots s_{l-p_0}, & k=p_0. \end{array} \right.  \] 
\end{lem} 
{\sc Proof.} We prove by induction over $k$. Assume $k+1<p_0$, then
\[ d^{k+1}(b)\stackrel{IH}{=} d(\Delta ^{p+k}s_{k+1}\cdots s_{l-k})=\Delta ^{p+k} \tau ^{p+k}(s_{l-k})s_{k+1}\cdots s_{l-k-1}. \]
Since $k+1<p_0$ we have $\inf (d^{k+1}(b))\le p+k+1<p+p_0=\inf _s(b)$, i.e. $d^{k+1}(b) \notin SSS(b)$. By Proposition \ref{Decyc}
we conclude that  $\tau ^{p+k}(s_{l-k})s_{k+1} =\Delta $. Hence $d^{k+1}(b)=\Delta ^{p+(k+1)}s_{(k+1)+1}\cdots s_{l-(k+1)}$. \\
Note that $d^k(b)$ may not be in $SSS(b)$. Indeed it is if and only if the pair $(\tau ^{p+k}(s_{l-k}), s_{k+1})$ is left-greedy. 
If it is not, then we obtain an SSS-element by one more decycling which decreases the supremum. \\
In the case $k+1=p_0$ we also may show by induction that 
\begin{eqnarray*}
 d^{p_0}(b) &=& d(d^{p_0-1}(b))\stackrel{IH}{=} d(\Delta ^{p+p_0-1}s_{p_0}\cdots s_{l-p_0+1})  \\
 &=& \Delta ^{p+p_0-1} \tau ^{p+p_0-1}(s_{l-p_0+1})s_{p_0}s_{p_0 +1}\cdots s_{l-p_0}. 
\end{eqnarray*}
We conclude that $\Delta \preceq \tau ^{p+p_0-1}(s_{l-p_0+1}) s_{p_0}$ - otherwise $\inf _s(b)=p+p_0-1$ by Propositions \ref{SSSrigid} and \ref{Decyc}.
Write  $\tau ^{p+p_0-1}(s_{l-p_0+1}) s_{p_0}=\Delta s_{p_0}''$ and we get the assertion. \\
Note that here are still two cases to consider. If $s_{p_0}''\ne 1$ then $d^{p_0}(b)\in SSS(b)$ according to Proposition \ref{Decyc}.
But if $s_{p_0}''=1$ then $d^{p_0}(b)$ is in $SSS(b)$ if and only if the pair $(\tau ^{p+p_0}(s_{l-p_0}), s_{p_0+1})$ is left-greedy.
Now if that pair is not left-greedy, then we may obtain an SSS-element by one more decycling which decreases the supremum. 
\quad $\Box $ 
 
\subsection{Super Summit Set}
 
Inside $SSS(b)$ (or $C(b)$) two elements can be conjugated by a sequence of conjugations by simple elements $s\in S$. Therefore, one may define
a conjugacy graph whose vertices are the elements of a conjugacy class (or of $SSS(b)$) and the edges are conjugations by simple elements $s\in S$.
In general, the SSS might be quite big. Certainly no polynomial bound (in $\ell (b)$) on the size (or even on the diameter) of the SSS graph is known. 
But in the case of Garside systems {\it of spindle type} the situation is fortunately simpler.
Proposition \ref{StructSSS} completely describes the structure of SSS graphs
in Garside systems of spindle type.

\begin{prop} \label{StructSSS} \hspace{-0.15cm}{\bf .}
Let $(G,S)$ with $S={\rm Div}(\Delta )$ be a Garside system of spindle type.
Let $b\in SSS(b) \subset G$, and let $r={\rm ord}  (b)$ and $p=\inf (b)$.  \\ 
(a) Define $m=l \cdot lcm(p,r)/p$ for $p>0$ and $m=l$ for $p=0$.
Then the SSS-graph of $b$ is a quotient graph of the graph  shown in Figure \ref{fig:graph}.
\begin{figure}
\begin{center}$ \xymatrix{
 & \ar[d]^c & \ar[d]^c & \ar[d]^c & \ar[d]^c & \\ 
\ar[r]^{\Delta } & b \ar[r]^{\Delta } \ar[d]^c & \tau (b) \ar[r]^{\Delta } \ar[d]^c  & \ar@{.}[r] \ar[d]^c & \tau ^{r-1}(b) \ar[r]^{\Delta } \ar[d]^c & \\ 
\ar[r]^{\Delta } & c(b) \ar[r]^{\Delta } \ar@{.}[d] & c (\tau (b)) \ar[r]^{\Delta } \ar@{.}[d] &  \ar@{.}[r] \ar@{.}[d] & c(\tau ^{r-1}(b)) \ar[r]^{\Delta } \ar@{.}[d] & \\ 
\ar[r]^{\Delta } & c^{l-1}(b) \ar[r]^{\Delta } \ar[d]^c  & c^{l-1}(\tau (b)) \ar[r]^{\Delta } \ar[d]^c  & \ar@{.}[r] \ar[d]^c & c^{l-1}(\tau ^{r-1}(b)) \ar[r]^{\Delta } \ar[d]^c &  \\
\ar[r]^{\Delta } & \tau ^{-p}(b) \ar[r]^{\Delta } \ar@{.}[d]  & \tau ^{1-p}(b) \ar[r]^{\Delta } \ar@{.}[d]  & \ar@{.}[r] \ar@{.}[d] & \tau ^{r-p-1}(b) \ar[r]^{\Delta } \ar@{.}[d] &  \\ 
 & & & & &
} $\end{center} \caption{SSS-graph}
\label{fig:graph}\end{figure}
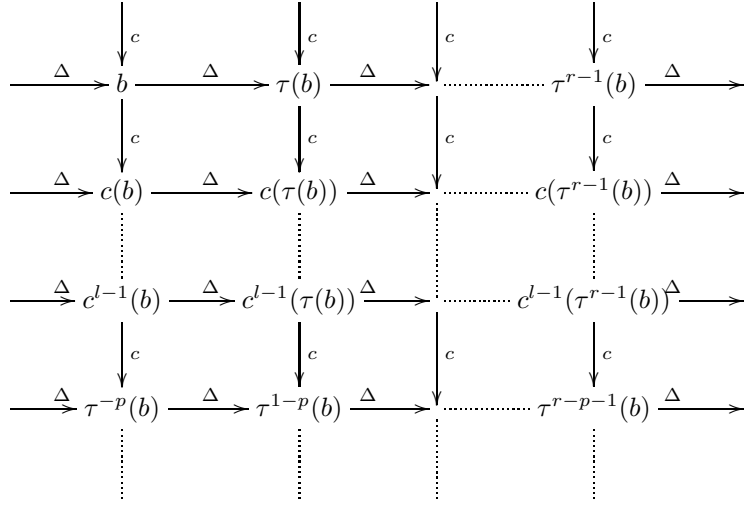
Note that we only draw cyclings $\xymatrix{\alpha \ar[r]^c &\beta }$, but this always implies also a decycling $\xymatrix{\alpha & \beta \ar[l]_d }$.
Also the graph is to be viewed as drawn on a torus.  \\
(b) We have the following bound on the size of the SSS: $|SSS(b)| \le l \cdot {\rm ord} (b)$ with $l=cl(b)$.
\end{prop}

{\sc Proof.} (a) First we show that in Garside systems of spindle type the conjugacy graph of the SSS may be obtained by repeated cyclings (conjugation by $\tau ^{-p}(s_1)$)
and $\tau $-automorphisms (conjugation by $\Delta $) only. It is well known that positive simple conjugations suffice. For $s\in S$ and $b\stackrel{LNF}{=}
\Delta ^ps_1\cdots s_l \in SSS(b)$, consider the conjugate $b^s=\Delta ^{p-1} \partial ^{2p-1} (s) s_1 \cdots s_l\cdot s$. 
We want to determine for which simple elements $s$, except for $1$ and $\Delta $, $b^s$ also belongs to $SSS(b)$. 
For $b^s\in SSS(b)$, we have $\inf (b^s)=p$. Therefore $\Delta \preceq \partial ^{2p-1} (s) s_1 \cdots s_l\cdot s$.
we conclude that either (a1) $\Delta \preceq \partial ^{2p-1} (s) s_1$ or (a2) $\Delta \preceq s_l s$.  \\
(a1) If $\partial ^{2p-1} (s) s_1=\Delta s_1''$ with $s_1=s_1's_1''$, we get $s_1'=\partial (\partial ^{2p-1}(s))=\tau ^p(s)$ $\Leftrightarrow $ $s=\tau ^{-p}(s_1')$. Hence
\[ b^s=\Delta ^{p-1} \partial ^{2p-1}(\tau ^{-p}(s_1'))s_1\cdots s_l\tau ^{-p}(s_1')=\Delta ^ps_1''s_2\cdots s_l \tau ^{-p}(s_1'). \]
It remains to show that $s_1''=1$, i.e. $s_1=s_1'$. Then $b^s=b^{\tau ^{-p}(s_1)}=c(b)$ is cycling. \\
Now, assume $s_1''\ne 1$. Since $s_1'=\tau ^p(s)$ and we consider $s\ne 1$, $s_1'$ is a nontrivial left divisor of $s_1$. Hence the pair
$(s_l, \tau ^{-p}(s_1'))$ is left-greedy iff $(s_l, \tau ^{-p}(s_1))$ is a left-greedy pair. And this is the case since $b$ is rigid as an element inside
$SSS(b)$. Also all other consecutive pairs in the product $s_1''s_2\cdots s_l \tau ^{-p}(s_1')$ are left-greedy except for the first. Indeed, if
$(s_1'', s_2)$ were left-greedy, then $cl(b^s)=l+1$ in contradiction with $b^s\in SSS(b)$. $\partial (s_1'')\wedge s_2 \ne 1$ implies either 
$\partial (s_1'')\wedge s_2 =\partial (s_1'')$ or $\partial (s_1'')\wedge s_2 =s_2 $. In the first case we get $\partial (s_1'') \preceq s_2$, hence
$\Delta \preceq s_1'' s_2$ and $\inf (b^s)=p+1$ in contradiction to $b^s\in SSS(b)$. In the latter case we have $s_2 \preceq \partial (s_1'')$, i.e. 
$s_1''s_2 \in S$. Now, since $1\ne s_1'' \preceq \partial (s_1')$, we have $\partial (s_1') \wedge (s_1''s_2)\ne 1$. This again splits in two cases.
 $\partial (s_1') \wedge (s_1''s_2)$ is either  $\partial (s_1')$ or  $s_1''s_2$. In the first case we have $\partial (s_1') \preceq s_1''s_2$,
 hence $\Delta \preceq s_1s_2$ in contradiction to $\inf (b)=p$. In the latter case we get $s_1'' s_2 \preceq \partial (s_1')$, hence $s_1s_2 \preceq
\Delta $ in contradiction with $\sup s_1s_2=2$. Thus we have proven $s_1''=1$. \\
(a2) If $\Delta \preceq s_l s$ we write $s=s's''$ with $s'=\partial (s_l)$. Hence $s''\preceq \partial (s')=\tau (s_l)$ and we write $s_l=s_l''s_l'$
with $s''=\tau (s_l'')$. Thus
\[ \partial ^{2p-1}(s)=\tau ^{p-1}(\partial (s))=\tau ^{p-1}((s'')^{-1}\partial (s'))=\tau^{p-1}((s'')^{-1}\tau (s_l))=\tau ^p(s_l'), \]
and we obtain
\[ b^s=\Delta ^{p-1}\tau ^p(s_l')s_1\cdots s_{l-1}\Delta s''=\tau [\Delta ^p \tau ^p(s_l')s_1\cdots s_{l-1}s_l'']. \]
Now if $s_l'=1$ then $b^s=\tau (b)$, and if $s_l''=1$ then we get $b^s=\tau (d(b))$, i.e. decycling up to a $\tau $-automorphism.
Analogeously as in (a1) one may show that the case $s_1'\ne 1$ and $s_1'' \ne 1$ contradicts to $b, b^s$ lying inside $SSS(b)$.  \\
Thus we have shown that cycling, decycling and $\tau $-operations suffice to enumerate $SSS(b)$.
Recall that, for $r={\rm ord}(b)$ and $p>0$, we have  $c^{lcm (p,r)l/p}(b)=\tau ^{-lcm (p,r)}(b)=b$. Hence $c^{lcm(p,r)l/p -1}(b)=c^{-1}(b)=d(b)$ for $p>0$,
and $c^{l-1}=d(b)$ for $p=0$. We conclude that cyclings and $\tau $-operations suffice. Since decyclings are inverse cyclings, obviously also only decyclings 
and $\tau $-operations suffice. \\
The assertion on the structure of the SSS-graph as well as claim (b) are simple corollaries of this result and the identities $c^l(b)=\tau ^{-p}(b)$ and $\tau ^r(b)=b$. \quad $\Box $

\begin{theo} \label{CPRigidGarsSys} \hspace{-0.15cm}{\bf .}
Let $(G,S)$ with $S={\rm Div}(\Delta )$ be a  Garside system of spindle type. If $G$ is a linear group, then the conjugacy problem is solvable in logspace.  
\end{theo} 

{\sc Proof.} Let $X$ be the set of atoms of the Garside system of spindle type $(G,S)$. We construct an L-computable function $f: (X^{\pm 1})^2 \longrightarrow
\{ {\rm true/false} \}$ that on input $(w,w')$ decides whether $\overline{w}$ and $\overline{w}'$ are conjugated in $G$. 
$f$ is the composition of three L-computable functions $f_a, f_b$ and $f_c$.

\noindent 
(a) Denote $S_{prop}=S \setminus \{ 1, \Delta \}$. The function $f_a: (X^{\pm 1})^2 \longrightarrow (\mathbb{Z} \times S_{prop})^2$ returns  the left Garside normal forms 
of  $\overline{w}$ and $\overline{w}'$, i.e., $(p,(s_1, \ldots , s_l),p',(s_1', \ldots , s_{l'}'))$ such that $\Delta ^ps_1 \cdots s_l=\overline{w}$ and  $\Delta ^{p'}s_1' \cdots s_{l'}'=\overline{w}'$ for some $l,l' \in \mathbb{N}$.   

\noindent
(b) Given the LNF's of  $\overline{w}$ and $\overline{w}'$, the function $f_b: (\mathbb{Z} \times S_{prop})^2 \longrightarrow (\mathbb{Z} \times S_{prop})^2$ computes 
the LNF's of elements $b,b' \in G$ such that $b\in SSS(\overline{w})$ and $b'\in SSS(\overline{w}')$.
Consider, for example,  the left normal form for $\overline{w}$, and assume that $\inf _s(\overline{w})=p+p_0$.
Lemma \ref{DecycPow} implies then that $\tau ^p(s_{l-p_0+1} \cdots s_l)s_1 \cdots s_{p_0-1}s_{p_0}'=\Delta ^{p_0}$ for some $s_{p_0}' \preceq s_{p_0}=s_{p_0}'s_{p_0}''$. Furthermore (see proof of Lemma \ref{DecycPow}), if $s_{p_0}'' \ne 1$ then $\Delta ^{p+p_0} s_{p_0}'' s_{p_0+1}\cdots s_{l-p_0}$ provides a left normal decomposition of an element in $SSS(\overline{w})$, namely of $d^{p_0}(\overline{w})$. But if $s_{p_0}''=1$ then
 $d^{p_0}(\overline{w})=\Delta ^{p+p_0} s_{p_0+1}\cdots s_{l-p_0} \in SSS(\overline{w})$ if $(\tau ^{p+p_0}(s_{l-p_0}), s_{p_0 +1})$ is left-greedy,
 and  
$$d^{p_0 +1}(\overline{w})=\Delta ^{p+p_0} (\tau ^{p+p_0}(s_{l-p_0}) s_{p_0+1}) \cdot s_{p_0+2} \cdots s_{l-p_0 -1} \in SSS(\overline{w})$$
 if  it is not, i.e. $(\tau ^{p+p_0}(s_{l-p_0}) s_{p_0+1}) \in S$. Therefore, given the LNF of  $\overline{w}$, in order to compute the LNF of an element
 $b\in SSS(\overline{w})$ one may perform the following steps. \\
(b1) Check, for convenience whether $\overline{w}$ is rigid. If yes, return $\overline{w}$ - otherwise proceed with (b2). \\
(b2) Check for all $p_0=l \, {\rm div} \, 2, \ldots , 1$ and for all prefixes $s_{p_0}'\preceq s_{p_0}$ whether
$\tau ^p(s_{l-p_0+1} \cdots s_l)s_1 \cdots s_{p_0-1}s_{p_0}' \stackrel{?}{=}\Delta ^{p_0}$. This can be done by the given WP-oracle.
Note that there exists a unique positive word for
the proper simple elements $s_{p_0}$. When found compute $s_{p_0}''=(s_{p_0}')^{-1} s_{p_0}$.  \\
(b3) If $s_{p_0}'' \ne 1$ return $(p+p_0, s_{p_0}'', s_{p_0+1}, \ldots , s_{l-p_0})$. If $s_{p_0}'' =1$ then
check whether  $(\tau ^{p+p_0}(s_{l-p_0}), s_{p_0 +1})$ is left-greedy and return the LNF as explained above. 

\noindent
(c) The function $f_c: (\mathbb{Z} \times S_{prop})^2 \longrightarrow \{ {\rm true/false} \}$ checks whether $b' \in SSS(b)$.
In abuse of notation $(p,(s_1, \ldots , s_l),p',(s_1', \ldots , s_{l'}'))$ encodes now the LNF's of $b$ and $b'$, respectively.
For convenience, we may first compare $(p,l)$ with $(p',l')$ and return false if they differ.
If not we compare for all $j=0,1, \ldots , l-1$ and for all $i=0,1,\ldots , r-1$
the LNF's of $\Delta ^{-p}b'$ with the LNF of $\Delta ^{-p} c^j(\tau ^i(b))$, i.e we compare
$(s_1', \ldots , s_{l'}')$ with $(\tau ^i(s_{j+1}, \ldots , \tau ^i(s_l), \tau ^{i-p}(s_1), \ldots , \tau ^{i-p}(s_j))$.
If we found a match, we return true, otherwise false.  
\quad $\Box $  \\

\noindent{\bf Remark.} It is easy to modify the algorithm in order to show that we may also solve the corresponding witness problem, namely the {\it conjugacy search problem} for linear Garside groups of spindle type, in logspace.

\section{Normal forms related to HNN extensions}  \label{VI}

The commutator subgroup $[B_n,B_n]$ of the $n$-strand braid group is the kernel of the abelianizer map $\mathcal{A}$ which sends every Artin generator
$\sigma _i$ to $t=\sigma _1$, i.e., we have the following short exact sequence.
\[ [B_n,B_n] \longrightarrow B_n \stackrel{\mathcal{A}}{\longrightarrow } \langle t \rangle.  \]
Indeed, $\langle \sigma _1 \rangle \le B_n / [B_n, B_n]$ is obvious, and the commutator relations $\sigma _i\sigma _{i+1}=\sigma _{i+1}\sigma _i$ (for 
$i=1, \ldots , n-2$) imply that $\sigma _1=\sigma _2= \ldots \sigma _{n-1}$ holds in $B_n / [B_n, B_n]$. \\ 
Gorin and Lin computed the commutator subgroups $[B_n, B_n]$ for all $n$ \cite{\GL} (see also \cite{\BV}).

For $n=3$, $[B_3, B_3]\cong F_2$ is freely generated, e.g. by $a=[\sigma _2,\sigma _1^{-1}]=\sigma _1\sigma _2\sigma _1^{-2}$ 
and $b=[\sigma _1,\sigma _2]=\sigma _2\sigma _1^{-1}$.
The $t^{\pm 1}$-action (by conjugation) on $\langle a, b \rangle $ is given by \cite{\BV}:
\[ t^{-1}a^{\pm 1}t=b^{\pm 1}, \quad t^{-1}b^{\pm 1}t=(ba^{-1})^{\pm 1}, \quad
ta^{\pm 1}t^{-1}=(b^{-1}a)^{\pm 1}, \quad tb^{\pm 1}t^{-1}=a^{\pm 1}. \]
In particular, $B_3$ admits the following presentation as an HNN extension over $F_2$.
\[ B_3=\langle a,b, t \mid at=tb, \quad bt=tba^{-1} \rangle . \]
Now, given a word over $\{ a,b,t \} ^{\pm 1}$ representing a braid, we may bring all powers off $t$ to the left.

\begin{defi} \label{HNNB3} \hspace{-0.15cm}{\bf .}
For every 3-strand braid $\beta \in B_3$ there exists an unique representation $t^p \cdot V(a,b)$ where $p\in \mathbb{Z}$ and $V=V(a,b)$ is a reduced 
group word over $\{a,b \}$. We call this representation the {\it (left) HNN normal form}  of $\beta $.
\end{defi}

One can show that this normal form is L-computable. Since we already know $B_3$ has logspace normal form for any finite generating set, we omit the details for this specific normal form here. 

\subsection{HNN normal form for $B_4$} \label{6.2}

The 4-strand braid group $B_4$ admits a description as tower of HNN extensions of the 2-rank free group, namely \cite{\BV}
\[ F_2=\langle a,b \rangle \subset [B_4, B_4]=\langle a,b,t_1,t_2 \rangle \subset B_4=\langle a,b,t_1,t_2,t \rangle \]
where $t_1=\sigma _1\sigma _2\sigma _1^{-2}=[\sigma _2,\sigma _1^{-1}]$, $t_2=\sigma _2\sigma _1^{-1}=[\sigma _1, \sigma _2]$,
$a=\sigma _1\sigma _2\sigma _1^{-1}\sigma _3\sigma _2^{-1} \sigma _1^{-1}=[\sigma _2\sigma _1^{-1},\sigma _3]=[\sigma _3\sigma _1^{-1},\sigma _2^{-1}]$,
and $b=\sigma _3\sigma _1^{-1}=[\sigma _2^{-1}\sigma _1^{-1}\sigma _3\sigma _2\sigma _3^2\sigma _1^{-2}, \sigma _1^{-1}]$. \\
The commutator subgroup $[B_4, B_4]$ is a semidirect product of the free group $\langle a, b \rangle $ with the free group 
$\langle t_1, t_2 \rangle \cong [B_3, B_3]$. And $\langle a, b \rangle $ is the kernel of the homomorphism
$B_4 \longrightarrow B_3$ \cite{\GL} given by
\[ t \mapsto t, \quad t_1 \mapsto t_1, \quad t_2 \mapsto t_2, \quad a \mapsto 1, \quad b \mapsto 1.  \]
With respect to these {\it Gorin-Lin generators} the 4-strand braid group admits the following presentation \cite{\BV}.
\[ B_4=\langle a,b,t_1,t_2,t \mid t_1t=tt_2=t_2tt_1, [a,t_1]=[b,t]=1, bat_2=t_2b \rangle.  \]
From that presentation one may derive the $t_1^{\pm 1}$- and $t_2^{\pm 1}$-action on $\langle a, b \rangle $ and the $t^{\pm 1}$-action on $a, b, t_1$, and $t_2$.
For the explicit equations we refer to \cite{\BV}. The conjugation action implies that one may bring the $t^{\pm 1}$'s to the very left, and the
$t_1^{\pm 1}$'s and $t_2^{\pm 1}$'s to the left of the $a$ and $b$'s.
\begin{defi} \label{HNNB4} \hspace{-0.15cm}{\bf .}
For every 4-strand braid $\beta $ there exists an unique representation $t^p \cdot VW$ where $p\in \mathbb{Z}$ and $V=V(t_1,t_2)$, $W=W(a,b)$ are reduced 
group words over $\{t_1, t_2\}$, $\{a, b\}$, respectively. We call this representation the {\it (left) HNN normal form} of $\beta \in B_4$.
\end{defi}
 
In the following we show that the HNN normal form of a 4-strand braid is {\it not} an L-computable function.

\begin{lem} \label{t1act} \hspace{-0.15cm}{\bf .}
Consider the family of 4-strand braids in $\langle a,b \rangle \subset B_4$ given by the action of $t_1$-powers upon $b$, i.e. 
the family $(b^{t_1^m})_{m\in \mathbb{N}}$.
Let  $|\cdot |_A$ denote the word length of the unique reduced free group word over the alphabet $A=\{ a,b \}$.
Then we have
\[ |(b)^{t_1^m}|_A=\Phi ^{2m+1} -\Phi^{-(2m+1)} \quad  where \quad \Phi=(1+\sqrt{5})/2 \] 
denotes the golden ratio.
\end{lem}

{\sc Proof.} The $t_1^{\pm 1}$-action on $\langle a,b \rangle $ is given by
\[ (a^{\pm 1})^{t_1}=b^{\pm 1}, \,\, (b^{\pm 1})^{t_1}=(ba^{-1}b^2)^{\pm 1}, \,\, (a^{\pm 1})^{t_1^{-1}}=(a^2b^{-1}a)^{\pm 1}, \,\,
(b^{\pm 1})^{t_1^{-1}}=a^{\pm 1}. \]
We observe that a reduced word for $b^{t_1^m}$ is a (semigroup) word over $\{a^{-1},b \}$ only.
Indeed, under the $t_1$-action $b$ is mapped to $ba^{-1}b^2$ and $a^{-1}$ to $b^{-1}$ which is always cancelled with some letter $b$ since each
occurence of $a^{-1}$ is always between two letters $b$. 
Denote by ${\bar \alpha }_m$ and $\beta _m$ the number of occurences of letters $a^{-1}$ and $b$,
inside the unique reduced word representing $b^{t_1^m}$.
We conclude that the following (recurrence) equations hold: ${\bar \alpha }_m=\beta _{m-1}$ (for $m \ge 1$) and
\[ \beta _m=3\beta _{m-1} - {\bar \alpha }_{m-1}=3\beta _{m-1}=3\beta _{m-1} -\beta _{m-2} \quad \forall m \ge 2. \]
Since  $|b^{t_1^m}|={\bar \alpha }_m +\beta _m$, also $(|b^{t_1^m}|_A)_{m\ge 0}$ satisfies the recurrence relation $|b^{t_1^m}|_A=3|b^{t_1^{m-1}}|_A-|b^{t_1^{m-2}}|_A$ for $m \ge 2$ with 
initial values $|b|_A=1$, $|b^{t_1}|_A=4$. \\
Now, it is an easy exercise to show that $|(b)^{t_1^m}|_A=\Phi ^{2m+1} -\Phi^{-(2m+1)}$. \quad $\Box $  \\

\noindent{\bf Remark.} This sequence (starting with $1,4,11,29,76,199,521,1364, \ldots $) is also known as {\it Bisection of Lucas sequence} (see OEIS A002878). 
%and it has the generating function $(1+x)/(1-3x+x^2)$. %Similar sequences will appear in section \ref{ArtNF}.

\begin{theo} \label{HNNB4NF} \hspace{-0.15cm}{\bf .}
The HNN normal form in $t^pV(t_1,t_2)W(a,b)$ of a 4-strand braid is not logspace computable.
\end{theo}
{\sc Proof.} Indeed, we show a slightly stronger result, namely that there exists no logspace transducer that computes the HNN normal form $V(t_1, t_2)W(a,b)$
of any 4-strand braid inside the commutator subgroup $[B_4, B_4]$. Indeed, the family of braids $(b^{t_1^m})_{m\in \mathbb{N}}$ has, according to Lemma \ref{t1act}, HNN normal forms of exponential length (w.r.t. the alphabet $\{t_1, t_2, a, b\}^{\pm 1}$) in $m$. Therefore, the time complexity for printing the HNN normal form of $b^{t_1^m}$
is not in ${\rm P}=DTIME(poly(m))$. Since $DSPACE(\log (m))={\rm L} \subseteq {\rm P}$, we conclude that there exists no such logspace transducer. 
\quad $\Box $   \\

It may be shown that, for any braid $\beta \in B_4$, we can find in logspace $p\in \mathbb{Z}$ and a group word
$U=U(t_1, t_2, a,b)$ such that $\beta =t^p \overline{U}$.  The problem remains to find a unique form for $U$ that could be computed in logspace.  \\

\noindent{\bf Remark.} One may show that several other classical normal forms for the 3-strand braid group can be also
computed in logspace. In particular the following normal forms are L-computable.

\begin{itemize}
\item Alternating normal form of $B_3$ (see e.g. chapter VII in \cite{\DD}).
\item Rotating normal form of $B_3$ (see e.g. chapter VIII in \cite{\DD}).
\item $\sigma $-positive/negative geodesic normal form for nontrivial braids in $B_3$.
\item Artin's combed normal form for $P_3$.
\end{itemize}

Furthermore, one may show that the natural normal form which we get from embedding $B_3$ in $Aut(F_3)$ is not L-computable.
Details are given in an extended version of this article uploaded to {\tt arXiv} \cite{EK13}\footnote{In this earlier version we used another terminology. In particular, we called Garside groups of spindle type {\it rigid Garside groups}.}.

\bibliographystyle{plain}
\bibliography{refs-logGarsideSubmitIJAC.bib}

\def\cprime{$'$}
\begin{thebibliography}{10}

\bibitem{MR0019087}
E.~Artin.
\newblock Theory of braids.
\newblock {\em Ann. of Math. (2)}, 48:101--126, 1947.

\bibitem{MR2032983}
David Bessis.
\newblock The dual braid monoid.
\newblock {\em Ann. Sci. \'Ecole Norm. Sup. (4)}, 36(5):647--683, 2003.

\bibitem{MR1815219}
Stephen~J. Bigelow.
\newblock Braid groups are linear.
\newblock {\em J. Amer. Math. Soc.}, 14(2):471--486 (electronic), 2001.

\bibitem{MR1654165}
Joan Birman, Ki~Hyoung Ko, and Sang~Jin Lee.
\newblock A new approach to the word and conjugacy problems in the braid
  groups.
\newblock {\em Adv. Math.}, 139(2):322--353, 1998.

\bibitem{MR1870512}
Joan~S. Birman, Ki~Hyoung Ko, and Sang~Jin Lee.
\newblock The infimum, supremum, and geodesic length of a braid conjugacy
  class.
\newblock {\em Adv. Math.}, 164(1):41--56, 2001.

\bibitem{MR2202556}
Leonid Bokut and Andrei Vesnin.
\newblock Gr\"obner-{S}hirshov bases for some braid groups.
\newblock {\em J. Symbolic Comput.}, 41(3-4):357--371, 2006.

\bibitem{MR0132791}
A.~H. Clifford and G.~B. Preston.
\newblock {\em The algebraic theory of semigroups. {V}ol. {I}}.
\newblock Mathematical Surveys, No. 7. American Mathematical Society,
  Providence, R.I., 1961.

\bibitem{MR1942303}
Arjeh~M. Cohen and David~B. Wales.
\newblock Linearity of {A}rtin groups of finite type.
\newblock {\em Israel J. Math.}, 131:101--123, 2002.

\bibitem{MR2463428}
Patrick Dehornoy, Ivan Dynnikov, Dale Rolfsen, and Bert Wiest.
\newblock {\em Ordering braids}, volume 148 of {\em Mathematical Surveys and
  Monographs}.
\newblock American Mathematical Society, Providence, RI, 2008.

\bibitem{MR1710165}
Patrick Dehornoy and Luis Paris.
\newblock Gaussian groups and {G}arside groups, two generalisations of {A}rtin
  groups.
\newblock {\em Proc. London Math. Soc. (3)}, 79(3):569--604, 1999.

\bibitem{MR2987387}
Volker Diekert, Jonathan Kausch, and Markus Lohrey.
\newblock Logspace computations in {C}oxeter groups and graph groups.
\newblock In {\em Computational and combinatorial group theory and
  cryptography}, volume 582 of {\em Contemp. Math.}, pages 77--94. Amer. Math.
  Soc., Providence, RI, 2012.

\bibitem{MR2004479}
Fran{\c{c}}ois Digne.
\newblock On the linearity of {A}rtin braid groups.
\newblock {\em J. Algebra}, 268(1):39--57, 2003.

\bibitem{MR1315459}
Elsayed~A. El-Rifai and H.~R. Morton.
\newblock Algorithms for positive braids.
\newblock {\em Quart. J. Math. Oxford Ser. (2)}, 45(180):479--497, 1994.

\bibitem{MR3030521}
Murray Elder, Gillian Elston, and Gretchen Ostheimer.
\newblock On groups that have normal forms computable in logspace.
\newblock {\em J. Algebra}, 381:260--281, 2013.

\bibitem{EK13}
Murray Elder and Arkadius Kalka.
\newblock Logspace computations for rigid {G}arside groups.
\newblock Arxiv:1310.0933.

\bibitem{MR1161694}
David B.~A. Epstein, James~W. Cannon, Derek~F. Holt, Silvio V.~F. Levy,
  Michael~S. Paterson, and William~P. Thurston.
\newblock {\em Word processing in groups}.
\newblock Jones and Bartlett Publishers, Boston, MA, 1992.

\bibitem{Go07}
Eddy Godelle.
\newblock Parabolic subgroups of {G}arside groups.
\newblock {\em J. of Algebra}, 317:1--16, 2007.

\bibitem{MR0251712}
E.~A. Gorin and V.~Ja. Lin.
\newblock Algebraic equations with continuous coefficients, and certain
  questions of the algebraic theory of braids.
\newblock {\em Mat. Sb. (N.S.)}, 78 (120):579--610, 1969.

\bibitem{MR1465024}
Eun~Sook Kang, Ki~Hyoung Ko, and Sang~Jin Lee.
\newblock Band-generator presentation for the {$4$}-braid group.
\newblock {\em Topology Appl.}, 78(1-2):39--60, 1997.
\newblock Special issue on braid groups and related topics (Jerusalem, 1995).

\bibitem{MR1888796}
Daan Krammer.
\newblock Braid groups are linear.
\newblock {\em Ann. of Math. (2)}, 155(1):131--156, 2002.

\bibitem{Lee10}
Eon-Kyung Lee.
\newblock A positive presentation for the pure braid group.
\newblock {\em Journal of the Chungcheong Mathematical Society},
  23(3):555--561, 2010.

\bibitem{MR0445901}
Richard~J. Lipton and Yechezkel Zalcstein.
\newblock Word problems solvable in logspace.
\newblock {\em J. Assoc. Comput. Mach.}, 24(3):522--526, 1977.

\bibitem{MR2340901}
Markus Lohrey and Nicole Ondrusch.
\newblock Inverse monoids: decidability and complexity of algebraic questions.
\newblock {\em Inform. and Comput.}, 205(8):1212--1234, 2007.

\bibitem{MR0264062}
Wilhelm Magnus and Ada Peluso.
\newblock On a theorem of {V}. {I}. {A}rnol\cprime d.
\newblock {\em Comm. Pure Appl. Math.}, 22:683--692, 1969.

\bibitem{Pi00}
Matthieu Picantin.
\newblock {\em {Petites groupes gaussiens}}.
\newblock PhD thesis, Universit\'{e} de Caen, 2000.

\bibitem{MR1842395}
Matthieu Picantin.
\newblock The conjugacy problem in small {G}aussian groups.
\newblock {\em Comm. Algebra}, 29(3):1021--1039, 2001.

\bibitem{MR563704}
Hans-Ulrich Simon.
\newblock Word problems for groups and contextfree recognition.
\newblock In {\em Fundamentals of computation theory ({P}roc. {C}onf.
  {A}lgebraic, {A}rith. and {C}ategorical {M}ethods in {C}omput. {T}heory,
  {B}erlin/{W}endisch-{R}ietz, 1979)}, volume~2 of {\em Math. Res.}, pages
  417--422. Akademie-Verlag, Berlin, 1979.

\bibitem{Tsa13}
Valdemar Tsanov.
\newblock Triangle groups, automorphic forms, and torus knots.
\newblock {\em L'Enseignement Math\'ematique}, 59(3-4):73--113, 2013.

\bibitem{Vass}
Svetla Vassileva.
\newblock The conjugacy problem in wreath products is decidable in log-space.
\newblock April 2013. AMS Sectional Meeting, Boston College.

\end{thebibliography}

\noindent{\it E-mail addresses:} {\tt murrayelder@gmail.com}, 
                        {\tt arkadius.kalka@rub.de}

\end{document}